\documentclass{birkjour}
\usepackage{amsthm}
\usepackage{pxfonts}
\usepackage{yfonts}
\usepackage{dsfont}
\usepackage{graphicx}
\usepackage{relsize}
\usepackage{xcolor}
\usepackage{enumerate}
\usepackage{amssymb}
\usepackage{mathrsfs}

 \newtheorem{thm}{Theorem}[section]
 
 \newtheorem{lemma}[thm]{Lemma}
 
 \theoremstyle{definition}
 \newtheorem{definition}[thm]{Definition}
 \theoremstyle{remark}
 \newtheorem{remark}[thm]{Remark}
 
 \numberwithin{equation}{section}

\def\bel{\begin{equation}\label}
	\def\eeq{\end{equation}}
\def\ds{\displaystyle}

\def\R{\mathbb R}
\def\Z{\mathbb Z}
\def\L{{\bf L}}
\def\bar{\overline}
\def\supp{{\bf supp~}}
\def\M{{\bf M}}

\begin{document}

%
%
%
%
%
%
%
%
%

\title[Triebel-Lizorkin spaces in Dunkl setting]
 {Triebel-Lizorkin spaces in Dunkl setting}

\author[Chuhan Sun]{Chuhan Sun}

\address{School of Mathematical Sciences \\
Zhejiang University \\
Hangzhou, 310058, Zhejiang \\
China}

\email{sunchuhan@westlake.edu.cn}


\author{Zhiming Wang}
\address{Institute for Theoretical Sciences \\
Westlake University \\
Hangzhou, 310030, Zhejiang \\
China}
\email{wangzhiming@westlake.edu.cn}


\date{June 27, 2025}

\begin{abstract}
We establish Triebel-Lizorkin spaces in the Dunkl setting which are associated with finite reflection groups on the Euclidean space. The group structures induce two nonequivalent metrics: the Euclidean metric and the Dunkl metric. In this paper, the $\L^2$ space and the Dunkl-Calder\'on-Zygmund singular integral operator in the Dunkl setting play a fundamental role. The main tools used in this paper are as follows: (i) the Dunkl-Calder\'on-Zygmund singular integral operator and a new Calder\'on reproducing formula in $\L^2$ with the Triebel-Lizorkin space norms; (ii) new test functions in terms of the $\L^2$ functions and distributions; (iii) the Triebel-Lizorkin spaces in the Dunkl setting which are defined by the wavelet-type decomposition and the analogous atomic decomposition of the Hardy spaces.
\end{abstract}

\maketitle

\section{Introduction}
\setcounter{equation}{0}

It is well-known that the Triebel-Lizorkin spaces play an important role in modern harmonic analysis. In the 1970's, the Triebel-Lizorkin space $\dot{F}^{\alpha,q}_{p}$ on $\R^n$ was introduced by several mathematicians. Lizorkin \cite{Lizorkin1}-\cite{Lizorkin2} and Triebel \cite{Triebel} independently investigated the Triebel-Lizorkin space $\dot{F}^{\alpha,q}_{p}(\R^n)$, $\alpha\in\R$, $1<p<\infty$, $1<q<\infty$.  Peetre \cite{Peetre1}-\cite{Peetre3} extended the range of the admissible parameters $p$ and $q$ to values less than one. In \cite{F-J1}, applying the basic representation formula of the form $f=\sum_Q\left<f,\varphi_Q\right>\psi_Q$ for a distribution $f$ on $\R^n$, the Triebel-Lizorkin space was defined through the Littlewood-Paley theory. More precisely, let $\varphi$, $\psi$ be functions on $\R^n$ satisfying: 
\begin{equation}
	\varphi,\psi\in\mathscr{S}(\R^n),
\end{equation}
\begin{equation}
	\supp\hat{\varphi},\supp\hat{\psi}\subseteq\left\{\xi\in\R^n:\frac{1}{2}\leq|\xi|\leq2\right\},
\end{equation}
\begin{equation}
	|\hat{\varphi}(\xi)|, |\hat{\psi}(\xi)|\geq c>0~~~\textrm{if}~~~\frac{3}{5}\leq|\xi|\leq\frac{5}{3},
\end{equation}
and
\begin{equation}
	\sum_{k\in\Z}\bar{\hat{\varphi}(2^{k}\xi)}\hat{\psi}(2^k\xi)=1~~~\textrm{if}~~~\xi\neq0.
\end{equation}

Set $\varphi_k(x)=2^{kn}\varphi(2^kx)$ and $\psi_k(x)=2^{kn}\psi(2^kx)$, $k\in\Z$. 
The Triebel-Lizorkin space $\dot{F}^{\alpha,q}_p(\R^n)$, $\alpha\in\R$, $0<p<\infty$, $0<q\leq\infty$, was defined by the collection of all $f\in\mathscr{S}'/\mathscr{P}(\R^n)$ (tempered distributions modulo polynomials) such that
\begin{equation}
	\|f\|_{\dot{F}^{\alpha,q}_p}=\left\|\left(\sum_{k\in\Z}\left(2^{k\alpha}\left|\varphi_k\ast f\right|\right)^q\right)^{1\over q}\right\|_{\L^p(\R^n)}<\infty,
\end{equation}
and when $q=\infty$, $\|f\|_{\dot{F}^{\alpha,\infty}_p}=\left\|\sup\limits_{k\in\Z}2^{k\alpha}\left|\varphi_k\ast f(x)\right|~\right\|_{\L^p(\R^n)}$.

As a generalization of $\R^n$, the space of homogeneous type was introduced by Coifman and Weiss in \cite{CW}, which provides a natural setting for studying function spaces. The homogeneous Triebel-Lizorkin spaces on spaces of homogeneous type were studied in \cite{Han}-\cite{Han and Sawyer}. More precisely, the Triebel-Lizorkin spaces on homogeneous type were introduced by using the family of operators $\{\textbf{D}_k\}_{k\in\Z}$ where $\textbf{D}_k=\textbf{S}_k-\textbf{S}_{k-1}$ and $\{\textbf{S}_k\}_{k\in\Z}$ is an approximation to the identity. For $|\alpha|<\varepsilon$, $1<p,q<\infty$, the Triebel-Lizorkin space $\dot{F}^{\alpha,q}_p$ is the collection of all $f\in\left(\mathcal{M}^{(\beta,\gamma)}\right)'$ with $0<\beta,\gamma<\varepsilon$ such that
\begin{equation}
	\|f\|_{\dot{F}^{\alpha,q}_p}=\left\|\left\{\sum_{k\in\Z}\left(2^{k\alpha}\left|\textbf{D}_k(f)\right|\right)^q\right\}^{\frac{1}{q}}\right\|_{\L^p(\R^n)}<\infty
\end{equation}
where $\mathcal{M}^{(\beta,\gamma)}$ is the space of test functions. See \cite{Han}-\cite{Han and Sawyer} for more details.

In recent years, the Dunkl setting is more and more important. In \cite{H-L-L-W}, the authors introduced the Triebel-Lizorkin spaces $\dot{F}^{\alpha,q}_{p}$ in the Dunkl setting for $1<p,q<\infty$. In this paper, we extend the range of $p,q$ to all $0<p,q<\infty$. The key tool is the discrete Calder\'on reproducing formula derived from the Dunkl-Poisson kernel (see \cite{Hp}). We mention that the Dunkl-Poisson kernel involves two nonequivalent metrics: the Euclidean metric and the Dunkl metric. 

Now we recall the framework of the Dunkl setting, see \cite{DH1}, \cite{DH2}, and \cite{TX}. In $\R^n$, the reflection $\sigma_\alpha$ with respect to the hyperplane $\alpha_\perp$ orthogonal to a nonzero vector $\alpha$ is given by $$\sigma_\alpha(x)=x-2\frac{\left<x,\alpha\right>}{\|\alpha\|^2}\alpha.$$ A finite set $R\subset\R^n\setminus\{0\}$ is called a root system if $\sigma_\alpha(R)=R$ for every $\alpha\in R$. Let $R$ be a root system in $\R^n$ normalized so that $\left<\alpha,\alpha\right>=2$ for $\alpha\in R$ and $G$ the finite reflection group generated by the reflections $\sigma_\alpha(\alpha\in R)$, where $\sigma_\alpha (x)=x-\left<\alpha,x\right>\alpha$ for $x\in\R^n$. Corresponding to this reflection group, we denote by $\mathcal{O}(x)$ the $G$-orbit of a point $x\in\R^n$. There is a natural metric between two $G$-orbits $\mathcal{O}(x)$ and $\mathcal{O}(y)$, given by
\bel{1.7}
d(x,y):=\min_{\sigma\in G}\|x-\sigma(y)\|.
\eeq
Obviously, $d(x,y)\leq\|x-y\|$, $d(x,y)=d(y,x)$ and $d(x,y)\leq d(x,z)+d(z,y)$ for all $x,y,z\in\R^n$.

A multiplicity function $\kappa$ defined on $R$ (invariant under $G$) is fixed throughout this paper. Let 
\bel{}
d\omega(x)=\prod_{\alpha\in R}\left|\left<\alpha,x\right>\right|^{\kappa(\alpha)}dx
\eeq
be the associated measure in $\R^n$, see \cite{JA}, where, here and subsequently, $dx$ stands for the Lebesgue measure in $\R^n$. We denote by $N=n+\sum_{\alpha\in R}\kappa(\alpha)$ the homogeneous dimension of the system, here $N$ is also called the upper dimension. Observe that for $x\in\R^n$ and $r>0$, 
\bel{sim}
\omega(B(x,r))\sim r^n\prod_{\alpha\in R}\left(\left|\left<\alpha,x\right>\right|+r\right)^{\kappa(\alpha)},
\eeq 
hence $\inf_{x\in\R^n}\omega(B(x,1))\geq C>0$, and $\omega(B(x,r))\geq Cr^N$. According to (\ref{1.7}), we also have $\omega(B(x,r))\sim\omega(B(y,r))$ when $d(x,y)\sim r$, and $\omega(B(x,r))\leq\omega(B_d(x,r))\leq|G|\omega(B(x,r))$, where $B_d(x,r):=\left\{y\in\R^n:d(x,y)<r\right\}$.

Moreover, the measure $d\omega(x)$ satisfies
\bel{}
C^{-1}\left(\frac{r_2}{r_1}\right)^n\leq\frac{\omega(B(x,r_2))}{\omega(B(x,r_1))}\leq C\left(\frac{r_2}{r_1}\right)^N ~~~~~\textrm{for}~0<r_1<r_2.
\eeq
This implies that $d\omega(x)$ satisfies the doubling and reverse doubling properties, that is, there exists a constant $C>0$ such that for all $x\in\R^n$, $r>0$ and $\lambda\geq1$,
\bel{}
C^{-1}\lambda^n\omega(B(x,r))\leq\omega(B(x,\lambda r))\leq C \lambda^N \omega(B(x,r)).
\eeq

Next we consider the Dunkl (differential) operators $\textbf{T}_j$ defined by
\bel{}
\textbf{T}_jf(x)=\partial_jf(x)+\sum_{\alpha\in R} \frac{\kappa(\alpha)}{2}\left<\alpha,e_j\right>\frac{f(x)-f(\sigma_\alpha(x))}{\left<\alpha,x\right>},
\eeq
where $e_1$, $e_2$, $\dots$, $e_n$ are standard unit vectors of $\R^n$. The Dunkl Laplacian is then defined as $\Delta_D=\sum_{j=1}^n\textbf{T}_j^2$, which is equivalent to
\bel{}
\Delta_Df(x)=\Delta_{\R^n}f(x)+\sum_{\alpha\in R}\kappa(\alpha) \left(\frac{\partial_\alpha f(x)}{\left<\alpha,x\right>}-\frac{f(x)-f(\sigma_\alpha(x))}{\left<\alpha,x\right>^2}\right).
\eeq
Here $\Delta_{\R^n}$ is the standard Euclidean Laplacian.

The operator $\Delta_D$ is self-adjoint on $\L^2(\R^n,\omega)$, see \cite{JA}, and generates the heat semigroup
\bel{}
\textbf{H}_tf(x)=e^{t\Delta_D}f(x)=\int_{\R^n}H_t(x,y)f(y)d\omega(y),
\eeq
where the heat kernel $H_t(x,y)$ is a $C^\infty$ function for all $t>0$, $x,y\in\R^n$ and satisfies $H_t(x,y)=H_t(y,x)>0$ and $\int_{\R^n} H_t(x,y)d\omega(y)=1$.

The Poisson semigroup $\textbf{P}_t=e^{-t\sqrt{-\Delta_D}}$ is subordinated to the heat semigroup $\textbf{H}_t=e^{t\Delta_D}$ by
\bel{}
\textbf{P}_tf(x)=\pi^{-1/2}\int_0^\infty e^{-u}e^{\frac{t^2}{4u}\Delta_D}f(x)\frac{du}{\sqrt{u}}
\eeq
and correspondingly for their integral kernels
\bel{}
P_t(x,y)=\pi^{-1/2}\int_0^{2\pi}e^{-u}H_{\frac{t^2}{4u}}(x,y)\frac{du}{u}.
\eeq
Moreover, $u(x,t)=\textbf{P}_tf(x)$, so-called the Dunkl Poisson integral, solves the boundary value problem
\bel{}
\left\{
\begin{array}{ll}
	(\partial_t^2+\Delta_D)u(x,t)=0  \vspace{2mm} \\ 
	u(x,0)=f(x)
\end{array}
\right.
\eeq
in the half-space $\R^{n+1}_+$, see \cite{JA}.

Observe that $\{\textbf{P}_{2^{-k}}\}_{k\in\Z}$ is an approximation to the identity on $\L^2(\R^n,\omega)$, that is, for $f\in\L^2(\R^n,\omega)$
\bel{}
\lim_{k\rightarrow+\infty}\textbf{P}_{2^{-k}}(f)(x)=f(x),~~\lim_{k\rightarrow-\infty}\textbf{P}_{2^{-k}}(f)(x)=0.
\eeq

$\diamond$ Throughout, we denote the operator $\textbf{T}$ and its kernel $T(x,y)$ by the same letter with bold type text and plain text respectively, such that
$$
\textbf{T}(f)(x)=\int_{\R^n}T(x,y)f(y)d\omega(y).
$$

Now we consider the Triebel-Lizorkin spaces in the Dunkl setting. In \cite{Bui}, by using the Dunkl Laplacian $\Delta_D$, the author proved that the Triebel-Lizorkin spaces associated with the Dunkl Laplacian are identical to the Triebel-Lizorkin spaces defined in the space of homogeneous type $(\R^n,\|\cdot\|,\omega)$. 

In this paper, we establish the Triebel-Lizorkin spaces in the Dunkl setting in a different way. The critical difference is we use the Calder\'on reproducing formula derived from two different kernels. Our the new method is using a new Calder\'on reproducing formula in $\L^2$ (see \cite{Hp}) with the Triebel-Lizorkin space norms derived from the Dunkl-Poisson kernel. Set $$\textbf{D}_k=\textbf{P}_{2^{-k}}-\textbf{P}_{2^{-k-1}}$$ with the kernel $D_k(x,y)=P_{2^{-k}}(x,y)-P_{2^{-k-1}}(x,y)$. Then applying the Coifman's decomposition (see \cite{D-J-S}) of the identity on $\L^2(\R^n,\omega)$, we have
\bel{}
\textbf{I}=\left(\sum_{\ell\in\Z}\textbf{D}_\ell\right)\left(\sum_{k\in\Z}\textbf{D}_k\right)
=\sum_{|k-\ell|\leq M}\textbf{D}_\ell \textbf{D}_k+\sum_{|k-\ell|>M}\textbf{D}_\ell \textbf{D}_k
=\sum_{k\in\Z}\textbf{D}_k^M\textbf{D}_k+\textbf{R}_1
\eeq
where $M$ is a fixed constant, and
\begin{equation}
	\textbf{D}_k^M=\sum_{\ell:|k-\ell|\leq M}\textbf{D}_\ell, ~~~~~\textbf{R}_1=\sum_{|k-\ell|>M}\textbf{D}_\ell \textbf{D}_k.
\end{equation}

Therefore, for $f\in\L^2(\R^n,\omega)$ we have 
\bel{}
f(x)=\sum_{k\in\Z}\textbf{D}_k^M\textbf{D}_k(f)(x)+\textbf{R}_1(f)(x)=\textbf{T}_M(f)(x)+\textbf{R}_1(f)(x)+\textbf{R}_2(f)(x),
\eeq
where 
\bel{TM}
\textbf{T}_M(f)(x)=\sum_{k\in\Z}\sum_{Q\in Q_d^k}\omega(Q)D_{k}^M(x,x_Q)\textbf{D}_k(f)(x_Q),
\eeq
\bel{R_2}
\textbf{R}_2(f)(x)=\sum_{k\in\Z}\sum_{Q\in Q_d^k}\left(\int_QD_k^M(x,y)\textbf{D}_k(f)(y)d\omega(y)-\omega(Q)D_{k}^M(x,x_Q)\textbf{D}_k(f)(x_Q)\right),
\eeq
where $Q_d^k$ is the collection of all dyadic cubes $Q$ with the side length $2^{-M-k}$, $M$ is some fixed large integer, and $x_Q$ is any fixed point in the cube $Q$.

By showing the operator $\textbf{T}_M$ is invertible, the authors in \cite{Hp} proved the Calder\'on reproducing formula on $\L^2$ as following
\begin{thm}\label{L2}{\rm{\cite{Hp}}} 
	If $f\in\L^2(\R^n,\omega)$, then there exists a function $h\in\L^2(\R^n,\omega)$, such that $\|f\|_{2}\sim\|h\|_{2}$ and
	\bel{}
	f(x)=\sum_{k\in\Z}\sum_{Q\in Q_d^k}\omega(Q)D_k^M(x,x_Q)\textbf{D}_k(h)(x_Q)
	\eeq
	where $Q$ are dyadic cubes in $\R^n$, $Q_d^k$ is the collection of $Q$ with the side length $\ell(Q) = 2^{-k-M}$, and $x_Q$ are any fixed point in $Q$.
\end{thm}
Based on the above theorem, we introduce the following 
\begin{definition}\label{L-P q}
	Suppose that $f\in\L^2(\R^n,\omega)$, $\alpha\in\R$, $0<p<\infty$, $0<q<\infty$, the Littlewood-Paley $q$-function $\textbf{S}^\alpha_q(f)$ is defined by 
	\bel{}
	\textbf{S}^\alpha_q(f)(x):=\left\{\sum_{k\in\Z}\sum_{Q\in Q_d^k}\left(2^{k\alpha}\left|\textbf{D}_k(f)(x_Q)\right|\right)^q\chi_Q(x)\right\}^{1\over q},
	\eeq
	where $\chi_Q(x)$ is the characteristic function of the cube $Q$.
\end{definition}

To establish the Triebel-Lizorkin space in the Dunkl setting, we utilize the new Calder\'on reproducing formula in $\L^2$ with the Triebel-Lizorkin space norms derived from the Dunkl-Poisson kernel (see \cite{JA}). Next, we introduce the new test functions in terms of the $\L^2$ functions and distributions by the duality estimates, which is crucial for developing the Dunkl-Triebel-Lizorkin spaces. Then we consider the Dunkl-Triebel-Lizorkin spaces as the collection of some distributions on the test function, and define the spaces by the wavelet-type decomposition and the analogous atomic decomposition of the Hardy spaces. Finally,  we prove the Dunkl-Triebel-Lizorkin spaces defined this way are complete. 

\section{Formulation on the main results}
\setcounter{equation}{0}

The discrete Littlewood-Paley $q$-function in Definition \ref{L-P q} leads to introduce the Dunkl-Triebel-Lizorkin space norm for $f\in\L^2(\R^n,\omega)$ as follows:
\begin{definition}\label{T-L norm}
	For $f\in\L^2(\R^n,\omega)$, the Dunkl-Triebel-Lizorkin space norm of $f$ is defined as
	\bel{}
	\|f\|_{\dot{F}^{\alpha,q}_p}:=\|\textbf{S}^{\alpha}_q(f)\|_{p}
	\eeq
	for $\alpha\in\R$, $0<p<\infty$ and $0<q<\infty$.
\end{definition}

Here and subsequently, for $0<p<\infty$, we denote the norm of $f\in\L^p(\R^n,\omega)$ by
\[
\|f\|_p:=\|f\|_{\L^p(\R^n,\omega)}=\left\{\int_{\R^n}|f(x)|^pd\omega(x)\right\}^{1\over p}<\infty.
\]

Our first main result, the discrete Calder\'on reproducing formula for $f\in\L^2(\R^n,\omega)$ with respect to the Dunkl-Triebel-Lizorkin space norm, is given by the following
\begin{thm}\label{Thm 3}
	If $f\in\L^2(R^n,\omega)$ with $\|f\|_{\dot{F}^{\alpha,q}_p}<\infty$, for $|\alpha|<1$, $\max\left\{\frac{N}{N+1},\frac{N}{N+\alpha+1}\right\}<p<\infty$, $\max\left\{\frac{N}{N+1},\frac{N}{N+\alpha+1}\right\}<q<\infty$, where $N$ is the upper dimension, then there exists a function $h\in\L^2(\R^n,\omega)$, such that $\|f\|_2\sim\|h\|_2$, $\|f\|_{\dot{F}^{\alpha,q}_p}\sim\|h\|_{\dot{F}^{\alpha,q}_p}$ and
	\bel{reproducing formula in T-L}
	f(x)=\sum_{k\in\Z}\sum_{Q\in Q_d^k}\omega(Q)D_k^M(x,x_Q)\textbf{D}_k(h)(x_Q)
	\eeq
	where $Q_d^k$ is the collection of all dyadic cubes $Q$ with the side length $2^{-M-k}$, $M$ is some fixed large integer, and $x_Q$ is any fixed point in the cube $Q$.
	The series converges in $\L^2(\R^n,\omega)$ norm and the Dunkl-Triebel-Lizorkin space norm.
\end{thm}

Applying the above theorem, we provide the duality estimates which will be a key idea for developing the Dunkl-Triebel-Lizorkin space theory. Before we state the estimates, we first define some necessary space norms.

\begin{definition}\label{CMOapq}
	Suppose that $|\alpha|<1$, $0<p\leq 1$, $1<q<\infty$ and $f\in \L^2(\R^n, \omega)$. The norm of $f$ in the Dunkl-Carleson measure space $CMO^{\alpha, q}_p$ is defined by
	\bel{definition of CMOp}
	\|f\|_{CMO_{p}^{\alpha,q}} = \sup_P \left( \frac{1}{\omega(P)^{\frac{q}{p}-\frac{q}{q'}}} \sum_{Q\in P}\omega(Q)\left|2^{k\alpha}\textbf{D}_k^M(f)(x_Q)\right|^q\right)^{\frac{1}{q}}
	\eeq
	where $P$ runs over all dyadic cubes. 
\end{definition}

\begin{definition}\label{Fap infty}
	Suppose that $|\alpha|<1$, $1\leq p<\infty$ and $f\in\L^2(\R^n,\omega)$. The norm of $f\in\dot{F}^{\alpha,\infty}_{p}(\R^n,\omega)$ is defined by
	\bel{}
	\|f\|_{\dot{F}^{\alpha,\infty}_{p}}=\left\|{\sup_{{k\in\Z},{Q\in Q_d^k}}}2^{k\alpha}\left|\textbf{D}_k^M(f)(x_Q)\right|\chi_{Q}(x)\right\|_{p}.
	\eeq
\end{definition}

\begin{definition}\label{Fa infty infty}
	Suppose that $|\alpha|<1$ and $f\in\L^2(\R^n,\omega).$ The norm of $f\in\dot{F}^{\alpha,\infty}_{\infty}(\R^n,\omega)$ is defined by
	\bel{}
	\|f\|_{\dot{F}^{\alpha,\infty}_{\infty}}={\sup_{{k\in\Z},{Q\in Q_d^k}}}2^{k\alpha}\left|\textbf{D}_k^M(f)(x_Q)\right|.
	\eeq
\end{definition}

Now we return to state the duality estimates
\begin{thm}\label{dual prop}
	Suppose that $f,g\in\L^2(\R^n,\omega)$, $p'$, $q'$ are the conjugates of $1<p,q<\infty$, respectively.
	\begin{enumerate}[(A)]
		\item \label{1<p,q<infty}
		$1<p<\infty$, $1<q<\infty$, there exists a constant $C$ such that
		\bel{1<p,q<infty}
		\left|\left<f,g\right>\right|\leq C \|f\|_{\dot{F}^{\alpha,q}_p}~\|g\|_{\dot{F}^{-\alpha,q'}_{p'}}
		\eeq
		\item \label{0<p<1,1<q<infty}
		$\max\left\{\frac{N}{N+1},\frac{N}{N+\alpha+1}\right\}<p\leq1$, $1<q<\infty$, there exists a constant $C$ such that
		\bel{0<p<1,1<q<infty}
		\left|\left<f,g\right>\right|\leq C \|f\|_{\dot{F}^{\alpha,q}_p}~\|g\|_{CMO^{-\alpha,q'}_{p}}
		\eeq
		where $\|g\|_{CMO^{-\alpha,q'}_{p}}$ is defined as in Definition \ref{CMOapq}.
		\item \label{1<p<infty,0<q<1}
		$1<p<\infty$, $\max\left\{{N\over N+1}, {N\over N+\alpha+1}\right\}<q\leq 1$, there exists a constant $C$ such that
		\bel{1<p<infty,0<q<1}
		\left|\left<f,g\right>\right|\leq C \|f\|_{\dot{F}^{\alpha,q}_p}~\|g\|_{\dot{F}^{-\alpha,\infty}_{p'}}
		\eeq
		where $\|g\|_{\dot{F}^{-\alpha,\infty}_{p'}}$ is defined as in Definition \ref{Fap infty}.
		\item\label{0<p,q<1}
		$\max\left\{{N\over N+1}, {N\over N+\alpha+1}\right\}<p\leq1$, $\max\left\{{N\over N+1}, {N\over N+\alpha+1}\right\}<q\leq 1$, there exists a constant $C$ such that
		\bel{0<p,q<1}
		\left|\left<f,g\right>\right|\leq C \|f\|_{\dot{F}^{\alpha,q}_p}~\|g\|_{\dot{F}^{-\alpha+N({1\over p}-1),\infty}_{\infty}}
		\eeq
		where $\|g\|_{\dot{F}^{-\alpha+N({1\over p}-1),\infty}_{\infty}}$ is defined as in Definition \ref{Fa infty infty}.
	\end{enumerate}
\end{thm}

The above Theorem \ref{dual prop} means that for $1<p<\infty$, $1<q<\infty$, each function $f\in\L^2(\R^n,\omega)$ with $\|f\|_{\dot{F}^{\alpha,q}_p}<\infty$ can be considered as a continuous linear functional on $\L^2(\R^n,\omega)\cap\dot{F}^{-\alpha,{q'}}_{p'}(\R^n,\omega)$, the subspace of $g\in\L^2(\R^n,\omega)$ with the norm $\|g\|_{\dot{F}^{-\alpha,q'}_{p'}}<\infty$.

Therefore, one can consider $\L^2(\R^n,\omega)\cap\dot{F}^{-\alpha,q'}_{p'}(\R^n,\omega)$ as a new test function space and write $f\in \left(\L^2(\R^n,\omega)\cap\dot{F}^{\alpha,q'}_{p'}(\R^n,\omega)\right)'$, where $\left(\L^2(\R^n,\omega)\cap\dot{F}^{\alpha,q'}_{p'}(\R^n,\omega)\right)'$ is the distribution space. Other ranges of $p$, $q$ stated above have the same results. The following result describes an important property for such a distribution $f$. More precisely, we establish the following discrete Calder\'on reproducing formula in the distribution sense:
\begin{thm}\label{thm}
	For $|\alpha|<1$, $\max\left\{{N\over N+1}, {N\over N+\alpha+1}\right\}<p<\infty$, $\max\left\{{N\over N+1}, {N\over N+\alpha+1}\right\}<q<\infty$, suppose that $\{f_n\}_{n\in\Z}$ is a Cauchy sequence in $\L^2(\R^n,\omega)$ with $\|\textbf{S}^\alpha_q(f_n-f_m)\|_p\rightarrow 0$ as $m,n\rightarrow\infty$. Then there exists a distribution $f$ satisfies 
	\begin{enumerate}[(A)]
		\item \label{1<p,q<infty}
		For $1<p<\infty$, $1<q<\infty$, $f$ is a distribution in $ \left(\L^2(\R^n,\omega)\cap\dot{F}^{-\alpha,q'}_{p'}(\R^n,\omega)\right)'$ such that $$\lim\limits_{n\rightarrow\infty}\|f_n-f\|_{\dot{F}^{\alpha,q}_{p}}=0,$$
		and for each $g\in\L^2(\R^n,\omega)\cap\dot{F}^{-\alpha,q'}_{p'}(\R^n,\omega)$,
		$$\left<f,g\right>=\left<\sum_{k\in\Z}\sum_{Q\in Q_d^k}\omega(Q)D_k^M(\cdot,x_Q)\textbf{D}_k(h)(x_Q),g(\cdot)\right>$$
		where $h\in\left(\L^2(\R^n,\omega)\cap\dot{F}^{-\alpha,q'}_{p'}(\R^n,\omega)\right)'$ with $\|f\|_{\dot{F}^{\alpha,q}_p}\sim\|h\|_{\dot{F}^{\alpha,q}_p}$;  
		\item \label{0<p<1,1<q<infty}
		For $\max\left\{\frac{N}{N+1},\frac{N}{N+\alpha+1}\right\}<p\leq1$, $1<q<\infty$, $f$ is a distribution in $\left(\L^2(\R^n,\omega)\cap{CMO}^{-\alpha,q'}_{p}(\R^n,\omega)\right)'$, such that $$\lim\limits_{n\rightarrow\infty}\|f_n-f\|_{\dot{F}^{\alpha,q}_{p}}=0,$$
		and for each $g\in \L^2(\R^n,\omega)\cap{CMO}^{-\alpha,q'}_{p}(\R^n,\omega)$,
		$$\left<f,g\right>=\left<\sum_{k\in\Z}\sum_{Q\in Q_d^k}\omega(Q)D_k^M(\cdot,x_Q)\textbf{D}_k(h)(x_Q),g(\cdot)\right>$$
		where $h\in\left(\L^2(\R^n,\omega)\cap{CMO}^{-\alpha,q'}_{p}(\R^n,\omega)\right)'$ with $\|f\|_{\dot{F}^{\alpha,q}_p}\sim\|h\|_{\dot{F}^{\alpha,q}_p}$;
		\item \label{1<p<infty,0<q<1}
		For $1<p<\infty$, $\max\left\{{N\over N+1},{N\over N+\alpha+1}\right\}<q\leq1$, $f$ is a distribution in $\left(\L^2(\R^n,\omega)\cap\dot{F}^{-\alpha,\infty}_{p'}(\R^n,\omega)\right)'$ such that $$\lim\limits_{n\rightarrow\infty}\|f_n-f\|_{\dot{F}^{\alpha,q}_{p}}=0,$$
		and for each $g\in\L^2(\R^n,\omega)\cap\dot{F}^{-\alpha,\infty}_{p'}(\R^n,\omega)$,
		$$\left<f,g\right>=\left<\sum_{k\in\Z}\sum_{Q\in Q_d^k}\omega(Q)D_k^M(\cdot,x_Q)\textbf{D}_k(h)(x_Q),g(\cdot)\right>$$
		where $h\in\left(\L^2(\R^n,\omega)\cap\dot{F}^{-\alpha,\infty}_{p'}(\R^n,\omega)\right)'$ with $\|f\|_{\dot{F}^{\alpha,q}_p}\sim\|h\|_{\dot{F}^{\alpha,q}_p}$;
		\item \label{0<p,q<1}
		For $\max\left\{{N\over N+1},{N\over N+\alpha+1}\right\}<p\leq1$, $\max\left\{{N\over N+1}, {N\over N+\alpha+1}\right\}<q\leq1$, $f$ is a distribution \\ in $\left(\L^2(\R^n,\omega)\cap\dot{F}^{-\alpha+N(\frac{1}{p}-1),\infty}_{\infty}(\R^n,\omega)\right)'$ such that $$\lim\limits_{n\rightarrow\infty}\|f_n-f\|_{\dot{F}^{\alpha,q}_{p}}=0$$
		and for each $g\in\L^2(\R^n,\omega)\cap\dot{F}^{-\alpha+N(\frac{1}{p}-1),\infty}_{\infty}(\R^n,\omega)$,
		$$\left<f,g\right>=\left<\sum_{k\in\Z}\sum_{Q\in Q_d^k}\omega(Q)D_k^M(\cdot,x_Q)\textbf{D}_k(h)(x_Q),g(\cdot)\right>$$
		where $h\in\left(\L^2(\R^n,\omega)\cap\dot{F}^{-\alpha+N(\frac{1}{p}-1),\infty}_{\infty}(\R^n,\omega)\right)'$ with $\|f\|_{\dot{F}^{\alpha,q}_p}\sim\|h\|_{\dot{F}^{\alpha,q}_p}$.
	\end{enumerate}
\end{thm}

Then we define the Dunkl-Triebel-Lizorkin spaces by the following
\begin{definition}\label{definition}
	The Dunkl-Triebel-Linzorkin space $\dot{F}^{\alpha,q}_{p}(\R^n,\omega)$, $|\alpha|<1$, $\max\left\{{N\over N+1},{N\over N+\alpha+1}\right\}<p<\infty$, $\max\left\{{N\over N+1}, {N\over N+\alpha+1}\right\}<q<\infty$ is defined as follows:
	\bel{D-T space}
	\begin{array}{lr}\ds
		\dot{F}^{\alpha,q}_{p}=\left\{f:f(x)=\sum_{k\in\Z}\sum_{Q\in Q_d^k}\omega(Q)\lambda_QD_k^M(x,x_Q),
		~\textrm{with}~\left\|\left\{\sum_{k\in\Z}\sum_{Q\in Q_d^k}\left(2^{k\alpha}\left|\lambda_Q\right|\right)^q\chi_Q(x)\right\}^{1\over q}\right\|_{p}<\infty\right\}
	\end{array}
	\eeq
	where the series converges in the following distribution sense: 
	\begin{enumerate}[(A)]
		\item \label{1<p,q<infty}
		for $1<p<\infty$, $1<q<\infty$, $f$ converges in $\left(\L^2(\R^n,\omega)\cap\dot{F}^{-\alpha,q'}_{p'}(\R^n,\omega)\right)'$;  
		\item \label{0<p<1,1<q<infty}
		for $\max\left\{\frac{N}{N+1},\frac{N}{N+\alpha+1}\right\}<p\leq1$, $1<q<\infty$, $f$ converges in $\left(\L^2(\R^n,\omega)\cap{CMO}^{-\alpha,q'}_p(\R^n,\omega)\right)'$;
		\item \label{1<p<infty,0<q<1}
		for $1<p<\infty$, $\max\left\{{N\over N+1}, {N\over N+\alpha+1}\right\}<q\leq 1$, $f$ converges in $\left(\L^2(\R^n,\omega)\cap\dot{F}^{-\alpha,\infty}_{p'}(\R^n,\omega)\right)'$;
		\item \label{0<p,q<1}
		for $\max\left\{{N\over N+1},{N\over N+\alpha+1}\right\}<p\leq1$, $\max\left\{{N\over N+1},{N\over N+\alpha+1}\right\}<q\leq1$, $f$ converges \\ in $\left(\L^2(\R^n,\omega)\cap\dot{F}^{-\alpha+N(\frac{1}{p}-1),\infty}_{\infty}(\R^n,\omega)\right)'$.
	\end{enumerate}
	If $f\in\dot{F}^{\alpha,q}_{p}(\R^n,\omega)$, the norm of $f$ is defined by
	\bel{F norm}
	\|f\|_{\dot{F}^{\alpha,q}_p}:=\inf\left\{~\left\|~\left\{\sum_{k\in\Z}\sum_{Q\in Q_d^k}\left(2^{k\alpha}\left|\lambda_Q\right|\right)^q\chi_Q(x)\right\}^{1\over q}~\right\|_{p}~\right\}
	\eeq
	where the infimum is taken over all $f(x)=\sum_{k\in\Z}\sum_{Q\in Q_d^k}\omega(Q)\lambda_QD_k^M(x,x_Q)$.
\end{definition}

To clarify that the Dunkl-Triebel-Lizorkin space in Definition \ref{definition} is complete, we show the following
\begin{thm}\label{closure}
	For $|\alpha|<1$, $\max\left\{\frac{N}{N+1},\frac{N}{N+\alpha+1}\right\}<p<\infty$, $\max\left\{\frac{N}{N+1},\frac{N}{N+\alpha+1}\right\}<q<\infty$, 
	\bel{closure A}
	\dot{F}^{\alpha,q}_{p}(\R^n,\omega)=\bar{~\L^2(\R^n,\omega)\cap\dot{F}^{\alpha,q}_{p}(\R^n,\omega)~}.
	\eeq
	More precisely,
	\begin{enumerate}[(A)]
		\item \label{1<p,q<infty}
		for $1<p<\infty$, $1<q<\infty$, $\bar{~\L^2(\R^n,\omega)\cap\dot{F}^{\alpha,q}_{p}(\R^n,\omega)~}$ is the collection of all distributions $f\in\left(\L^2(\R^n,\omega)\cap\dot{F}^{-\alpha,q'}_{p'}(\R^n,\omega)\right)'$ such that there exists a sequence $\{f_n\}_{n=1}^{\infty}$ in $\L^2(\R^n,\omega)$ with $\|f_n-f_m\|_{\dot{F}^{\alpha,q}_{p}}\to0$ as $n,m\to\infty$. Moreover, $f_n$ converges to $f$ in $\left(\L^2(\R^n,\omega)\cap\dot{F}^{-\alpha,q'}_{p'}(\R^n,\omega)\right)'$;
		\item \label{0<p<1,1<q<infty}
		for $\max\left\{\frac{N}{N+1},\frac{N}{N+\alpha+1}\right\}<p\leq1$, $1<q<\infty$, $\bar{~\L^2(\R^n,\omega)\cap\dot{F}^{\alpha,q}_{p}(\R^n,\omega)~}$ is the collection of all distributions $f\in\left(\L^2(\R^n,\omega)\cap{CMO}^{-\alpha,q'}_{p}(\R^n,\omega)\right)'$ such that there exists a sequence $\{f_n\}_{n=1}^{\infty}$ in $\L^2(\R^n,\omega)$ with $\|f_n-f_m\|_{\dot{F}^{\alpha,q}_{p}}\to0$ as $n,m\to\infty$. Moreover, $f_n$ converges to $f$ \\ in $\left(\L^2(\R^n,\omega)\cap{CMO}^{-\alpha,q'}_{p}(\R^n,\omega)\right)'$;
		\item \label{1<p<infty,0<q<1}
		for $1<p<\infty$, $\max\left\{{N\over N+1},{N\over N+\alpha+1}\right\}<q\leq1$, $\bar{~\L^2(\R^n,\omega)\cap\dot{F}^{\alpha,q}_{p}(\R^n,\omega)~}$ is the collection of all distributions $f\in\left(\L^2(\R^n,\omega)\cap\dot{F}^{-\alpha,\infty}_{p'}(\R^n,\omega)\right)'$ such that there exists a sequence $\{f_n\}_{n=1}^{\infty}$ in $\L^2(\R^n,\omega)$ with $\|f_n-f_m\|_{\dot{F}^{\alpha,q}_{p}}\to0$ as $n,m\to\infty$. Moreover, $f_n$ converges to $f$ in $\left(\L^2(\R^n,\omega)\cap\dot{F}^{-\alpha,\infty}_{p'}(\R^n,\omega)\right)'$;
		\item \label{0<p,q<1}
		for $\max\left\{{N\over N+1},{N\over N+\alpha+1}\right\}<p\leq1$, $\max\left\{{N\over N+1}, {N\over N+\alpha+1}\right\}<q\leq1$, $\bar{~\L^2(\R^n,\omega)\cap\dot{F}^{\alpha,q}_{p}(\R^n,\omega)~}$ is the collection of all distributions $f\in\left(\L^2(\R^n,\omega)\cap\dot{F}^{-\alpha+N(\frac{1}{p}-1),\infty}_{\infty}(\R^n,\omega)\right)'$ such that there exists a sequence $\{f_n\}_{n=1}^{\infty}$ in $\L^2(\R^n,\omega)$ with $\|f_n-f_m\|_{\dot{F}^{\alpha,q}_{p}}\to0$ as $n,m\to\infty$. Moreover, $f_n$ converges to $f$ in $\left(\L^2(\R^n,\omega)\cap\dot{F}^{-\alpha+N(\frac{1}{p}-1),\infty}_{\infty}(\R^n,\omega)\right)'$. 
	\end{enumerate}
\end{thm}

The paper is organized as follows. In the next section, we prove Theorem \ref{Thm 3}, which is the Calder\'on reproducing formula in $\L^2\cap\dot{F}^{\alpha,q}_p$. The main tools are orthogonal estimates in the Dunkl setting. In section 4, we demonstrate the Theorem \ref{dual prop}, the duality estimates which lead a way for developing the Dunkl-Triebel-Lizorkin space theory. In the last section, we define the Triebel-Lizorkin space in the Dunkl setting in Definition \ref{definition} and show the spaces is complete by Theorem \ref{closure}.

\section{Calder\'on reproducing formula in $\L^2(\R^n,\omega)\cap\dot{F}^{\alpha,q}_p(\R^n,\omega)$}
\setcounter{equation}{0}

As mentioned before, the Dunkl-Calder\'on-Zygmund operator theory plays a crucial role. To prove Theorem \ref{Thm 3}, we recall the Dunkl-Calder\'on-Zygmund singular integral operator and almost orthogonality estimates in the Dunkl setting. See \cite{Hp} for more details.

Let $C_0^\eta(\R^n)$, $\eta>0$, denote the space of continuous functions $f$ with compact support and
\bel{}
\|f\|_{\eta}:=\sup_{x\neq y}\frac{|f(x)-f(y)|}{\|x-y\|^{\eta}}<\infty.
\eeq

\begin{definition}
	An operator $\textbf{T}:C_0^\eta(\R^n)\rightarrow(C_0^\eta(\R^n))'$ with $\eta>0$, is said to be a Dunkl-Calder\'on-Zygmund singular integral operator if $K(x,y)$, the kernel of $\textbf{T}$, satisfies the following estimates: for some $0<\varepsilon\leq1$,
	\bel{2.1 1}
	\left|K(x,y)\right|\leq \frac{C}{\omega(B(x,d(x,y)))}\left(\frac{d(x,y)}{\|x-y\|}\right)^\varepsilon
	\eeq
	for all $x\neq y$;
	\bel{2.1 2}
	\left|K(x,y)-K(x',y)\right|\leq \left(\frac{\|x-x'\|}{\|x-y\|}\right)^\varepsilon\frac{C}{\omega(B(x,d(x,y)))}
	\eeq
	for $\|x-x'\|\leq d(x,y)/2$;
	\bel{2.1 3}
	\left|K(x,y)-K(x,y')\right|\leq \left(\frac{\|y-y'\|}{\|x-y\|}\right)^\varepsilon\frac{C}{\omega(B(x,d(x,y)))}
	\eeq
	for $\|y-y'\|\leq d(x,y)/2$.
\end{definition}

A Dunkl-Calder\'on-Zygmund singular integral operator is said to be the Dunkl-Calder\'on-Zygmund operator if it extends a bounded operator on $\L^2(\R^n)$. Suppose that $\textbf{T}$ is the Dunkl-Calder\'on-Zygmund operator. We denote 
\bel{}
\|\textbf{T}\|_{dcz}=\|\textbf{T}\|_{2,2}+\|K\|_{dcz}
\eeq
as the Dunkl-Calder\'on-Zygmund operator norm, where $\|K\|_{dcz}$ is the minimum of the constants in (\ref{2.1 1})-(\ref{2.1 3}).

The following almost orthogonality estimates are important tools in the proof. See \cite{Hp} for more details. Let $\{\textbf{S}_k\}_{k\in\Z}$ be  approximations to the identity and set $\textbf{D}_k:=\textbf{S}_k-\textbf{S}_{k-1}$, then

\begin{lemma}\label{almost orthogonality 1}{\rm\cite{Hp}}
	For $k,j\in\Z$, $\varepsilon>0$, $\gamma$, $\varepsilon'\in(0,\varepsilon)$, 
	\bel{}
	\left|\left(D_kD_j\right)(x,y)\right|\leq C~2^{-|k-j|\varepsilon'}\frac{1}{V(x,y,2^{-k\vee-j}+d(x,y))}\left(\frac{2^{-k\vee-j}}{2^{-k\vee-j}+d(x,y)}\right)^\gamma,
	\eeq
	where $a\wedge b=\min\left\{a,b\right\}$, $a\vee b=\max\left\{a,b\right\}$, $V(x,y,r):=\max\left\{\omega(B(x,r)),\omega(B(y,r))\right\}$.
\end{lemma}

\begin{lemma}\label{almost orthogonality 2}{\rm \cite{Hp}}
	Let $\textbf{T}$ be a Dunkl-Calder\'on-Zygmund singular integral operator satisfying $\textbf{T}(1)=\textbf{T}^*(1)=0$ and $\textbf{T}$ is bounded in $\L^2(\R^n,\omega)$. Then
	\bel{}
	\begin{array}{lr}\ds
		~~~~\left|\int_{\R^n}\int_{\R^n}D_k(x,u)K(u,v)D_j(v,y)d\omega(u)d\omega(v)\right|
		\\\\ \ds 
		\leq C~2^{-|k-j|\varepsilon'} \|\textbf{T}\|_{dcz} \frac{1}{V(x,y,2^{-k\vee-j}+d(x,y))} \left(\frac{2^{-k\vee-j}}{2^{-k\vee-j}+d(x,y)}\right)^\gamma,
	\end{array}
	\eeq
	where  $\gamma$, $\varepsilon'\in(0,\varepsilon)$, $\varepsilon$ is the regularity exponent of the kernel of $\textbf{T}$ given in (\ref{2.1 2}) and (\ref{2.1 3}).
\end{lemma}

\textbf{Proof of Theorem \ref{Thm 3}:}
First, as mentioned before, we decompose the identity operator on $\L^2(\R^n,\omega)$ by $\textbf{I}=\textbf{T}_M+\textbf{R}_1+\textbf{R}_2$. To prove that $\textbf{T}_M$ is invertible and $(\textbf{T}_M)^{-1}$, the inverse of $\textbf{T}_M$, is bounded on $\L^2\cap\dot{F}^{\alpha,q}_p$, we need to estimate $\textbf{R}_1$ and $\textbf{R}_2$ on $\L^2\cap\dot{F}^{\alpha,q}_p$ and show that the norm of $\textbf{R}_1$ and $\textbf{R}_2$ on $\L^2\cap\dot{F}^{\alpha,q}_p$ are less than $1$. To this end, we consider the Dunkl setting, $(\R^n,\|\cdot\|,\omega)$, as a space of homogeneous type in the sense of Coifman and Wiess. The discrete Calder\'on reproducing formula in the space of homogeneous type is given by the following (see \cite{H-M-Y})
\begin{thm}\label{Thm 2}{\rm{\cite{H-M-Y}}}
	Let $\{\textbf{S}_k\}_{k\in\Z}$ be a Coifman's approximations to the identity and set $\textbf{E}_k:=\textbf{S}_k-\textbf{S}_{k-1}$. Then there exists a family of operators $\{\widetilde{\textbf{E}}_k\}_{k\in\Z}$ such that for any fixed $x_Q\in Q$ with $k\in\Z$ and $Q\in Q_{cw}^k$ are dyadic cubes with the side length $2^{-k-M_0}$,
	\bel{cwfE}
	f(x)=\sum_{k\in\Z}\sum_{Q\in Q_{cw}^k}\omega(Q)\widetilde{E}_k(x,x_Q)\textbf{E}_k(f)(x_Q),
	\eeq
	where the series converge in $\L^p(\R^n,\omega)$, $\mathcal{M}(\beta,\gamma,r,x_0)$ and $(\mathcal{M}(\beta,\gamma,r,x_0))'$, the dual space of $\mathcal{M}(\beta,\gamma,r,x_0)$. 
\end{thm}
Recall the Littlewood-Paley theory and the Triebel-Lizorkin spaces on space of homogeneous type $(\R^n, \|\cdot\|,\omega)$ in the sense of Coifman and Weiss. The discrete Calder\'on reproducing formula in Theorem \ref{Thm 2} leads the following discrete $q$-function on space of homogeneous type $(\R^n,\|\cdot\|,\omega)$:
\begin{definition}
	Suppose that $f\in(\mathcal{M}(\beta,\gamma,r,x_0))'$, $\alpha\in\R$, define the Littlewood-Paley $q$-function $\textbf{S}^\alpha_{q,cw}(f)$ for the space of homogeneous type $(\R^n,\|\cdot\|,\omega)$ as
	\bel{}
	\textbf{S}^\alpha_{q,cw}(f)(x)=\left\{\sum_{k\in\Z}\sum_{Q\in Q_{cw}^k}\left(2^{k\alpha}\left|\textbf{E}_k(f)(x_Q)\right|\right)^q\chi_Q(x)\right\}^{1\over q}~~~~\textrm{for}~0<q<\infty,
	\eeq
	where $\textbf{E}_k$, $Q_{cw}^k$ are given in Theorem \ref{Thm 2}.
\end{definition}

\begin{remark}
	It is well known that $\R^n$ together with the Euclidean metric and the Dunkl measure is space of homogeneous type in the sense of Coifman and Weiss. See \cite{H-M-Y} for more details.   
\end{remark}
Now we return to the estimates of $\textbf{R}_1$ and $\textbf{R}_2$. Note that $\textbf{R}_1$ and $\textbf{R}_2$ are Dunkl-Calder\'on-Zygmund operators with $\|\textbf{R}_1\|_{dcz}+\|\textbf{R}_2\|_{dcz}\leq C~2^{-M\delta}$, $\delta>0$. The boundedness of $\textbf{R}_1$ and $\textbf{R}_2$ on $\L^2$ follows from the Cotlar-Stein Lemma. Moreover, $\|\textbf{R}_1 + \textbf{R}_2\|_{2,2} < 1$. See \cite{Hp} for more details.

To estimate $\textbf{R}_1$ and $\textbf{R}_2$ on $\dot{F}^{\alpha,q}_p(\R^n,\omega)$, we establish the following estimates
\bel{}
\|\textbf{S}_q^\alpha(\textbf{R}_1(f))\|_p \leq C~\|\textbf{S}_{q,cw}^\alpha(\textbf{R}_1(f))\|_p\leq C~\|\textbf{R}_1\|_{dcz}\|\textbf{S}^\alpha_{q,cw}(f)\|_p\leq C~\|\textbf{R}_1\|_{dcz}\|\textbf{S}^\alpha_q(f)\|_p
\eeq 
and the similar estimates also hold for $\textbf{R}_2$. Now we show the above estimates for $\textbf{R}_1$ by the following steps:

{\bfseries Step 1: } $\|\textbf{S}_q^\alpha(\textbf{R}_1(f))\|_p \leq C~\|\textbf{S}_{q,cw}^\alpha(\textbf{R}_1(f))\|_p$.

Indeed we only need to show that for each $f\in \L^2(\R^n,\omega)$,
\bel{}
\|\textbf{S}_q^\alpha(f)\|_p \leq C~\|\textbf{S}_{q,cw}^\alpha(f)\|_p
\eeq
since $\textbf{R}_1$ is bounded on $\L^2$. By the discrete Calder\'on reproducing formula of $f\in\L^2$ given in Theorem \ref{Thm 2} we have
\bel{}
\begin{array}{lr}\ds
	\textbf{S}_q^\alpha(f)(x) = \left\{ \sum_{k'\in\Z} \sum_{Q'\in Q_d^{k'}} \left( 2^{k'\alpha} \left| \textbf{D}_{k'} \left( \sum_{k\in\Z} \sum_{Q\in Q_{cw}^k} \omega(Q) \widetilde{E}_k(\cdot,x_Q) \textbf{E}_k(f)(x_Q) \right) (x_{Q'}) \right|  \right)^q \chi_{Q'}(x) \right\}^{\frac{1}{q}}
	\\\\ \ds~~~~~~~~~~~~~~
	= \left\{ \sum_{k'\in\Z} \sum_{Q'\in Q_d^{k'}} \left(  \left| \sum_{k\in\Z} \sum_{Q\in Q_{cw}^k} 2^{k'\alpha} \omega(Q) (D_{k'} \widetilde{E}_k)(x_{Q'},x_Q) \textbf{E}_k(f)(x_Q) \right|  \right)^q \chi_{Q'}(x) \right\}^{\frac{1}{q}}.
\end{array}
\eeq

By the almost orthogonal estimate given in the Lemma \ref{almost orthogonality 1}, we choose $|\alpha|<\varepsilon<1$ such that
\bel{4.5}
\begin{array}{lr}\ds~~~
	\left|(D_{k'} \widetilde{E}_k)(x_{Q'},x_Q)\right|\chi_{Q'}(x)
	\\\\ \ds
	\leq C 2^{-|k'-k|\varepsilon} \frac{1}{V(x_{Q'},x_Q,2^{-k'\vee-k}+d(x_{Q'},x_Q))} \left(\frac{2^{-k'\vee-k}}{2^{-k'\vee-k}+d(x_{Q'},x_Q)}\right)^{\varepsilon} \chi_{Q'}(x)
	\\\\ \ds
	\leq C 2^{-|k'-k|\varepsilon} \frac{1}{\omega(B(x_Q,2^{-k'\vee-k}+d(x,x_Q)))} \left(\frac{2^{-k'\vee-k}}{2^{-k'\vee-k}+d(x,x_Q)}\right)^{\varepsilon} \chi_{Q'}(x).
\end{array}
\eeq
Since $d(x,y)=\min_{\sigma\in G}\|\sigma(x)-y\|$, then
\bel{4.6}
\begin{array}{lr}\ds~~~~
	\left|(D_{k'} \widetilde{E}_k)(x_{Q'},x_Q)\right|\chi_{Q'}(x)
	\\\\ \ds
	\leq C\sum_{\sigma\in G} 2^{-|k'-k|\varepsilon} \frac{1}{\omega(B(x_Q,2^{-k'\vee-k}+\|\sigma(x)-x_Q\|))} \left(\frac{2^{-k'\vee-k}}{2^{-k'\vee-k}+\|\sigma(x)-x_Q\|}\right)^{\varepsilon} \chi_{Q'}(x)
	\\\\ \ds
	\leq C\sum_{\sigma\in G} 2^{-|k'-k|\varepsilon} \frac{1}{\omega(B(\sigma(x),2^{-k'\vee-k}+\|\sigma(x)-x_Q\|))} \left(\frac{2^{-k'\vee-k}}{2^{-k'\vee-k}+\|\sigma(x)-x_Q\|}\right)^{\varepsilon} \chi_{Q'}(x).
\end{array}
\eeq
Let $\theta$ satisfies that $\max\left\{\frac{N}{N+\varepsilon},\frac{N}{N+\alpha+\varepsilon}\right\}<\theta<\min\{p, q, 1\}$, then
\bel{}
\begin{array}{lr}\ds~~~~~
	\left|\sum_{Q\in Q_{cw}^k} 2^{k'\alpha} \omega(Q) (D_{k'} \widetilde{E}_k)(x_{Q'},x_Q) \textbf{E}_k(f)(x_Q)\right|\chi_{Q'}(x)
	\\\\ \ds
	\leq C\sum_{Q\in Q_{cw}^k} \sum_{\sigma\in G} 2^{k'\alpha} 2^{-|k'-k|\varepsilon} \omega(Q) \frac{1}{\omega(B(\sigma(x),2^{-k'\vee-k}+\|\sigma(x)-x_Q\|))} \left(\frac{2^{-k'\vee-k}}{2^{-k'\vee-k}+\|\sigma(x)-x_Q\|}\right)^{\varepsilon} \left|\textbf{E}_kf(x_Q)\right|\chi_{Q'}(x)
	\\\\ \ds
	\leq C\sum_{\sigma\in G}\left\{  \sum_{Q\in Q_{cw}^k}  2^{-|k'-k|\varepsilon\theta} 2^{(k'-k)\alpha\theta}  \omega(Q)^\theta \frac{1}{\omega(B(\sigma(x),2^{-k'\vee-k}+\|\sigma(x)-x_Q\|))^\theta}\right.
	\\\\ \ds~~~~~ 
	\left.\times~\left(\frac{2^{-k'\vee-k}}{2^{-k'\vee-k}+\|\sigma(x)-x_Q\|}\right)^{\theta\varepsilon} \left|2^{k\alpha}\textbf{E}_kf(x_Q)\right|^\theta\right\}^{1\over \theta}\chi_{Q'}(x).
\end{array}
\eeq

Denote $c_Q$ as the center of $Q$. Set $$A_0=\left\{Q\in Q_{cw}^k:\|c_Q-\sigma(x)\|\leq 2^{-k\vee-k'}\right\}$$ and $$A_\ell=\left\{Q\in Q_{cw}^k: 2^{\ell-1+(-k\vee-k')}<\|c_Q-\sigma(x)\|\leq2^{\ell+(-k\vee-k')}\right\}$$ for $\ell\in\mathbb{N}$.

For $Q\in Q_{cw}^k$, we have 
\bel{}
\omega(Q)\chi_Q(z)\sim\omega(B(z,2^{-k-M_0}))\chi_Q(z)\sim\omega(B(\sigma(z),2^{-k-M_0}))\chi_Q(z),~{\textrm{for}}~\sigma\in G,
\eeq
and
\bel{estimate size of omega}
\omega(B(\sigma(z),2^{-k\vee-k'})) \leq C~2^{((-k\vee-k')-(-k-M_0))N}\omega(B(\sigma(z),2^{-k-M_0})).
\eeq
Hence
\bel{4.11}
\begin{array}{lr}\ds
	~~~\sum_{Q\in Q_{cw}^k} 2^{-|k'-k|\varepsilon\theta} 2^{(k'-k)\alpha\theta} \omega(Q)^\theta \frac{1}{\omega(B(\sigma(x),2^{-k'\vee-k}+\|\sigma(x)-x_Q\|))^\theta} \left(\frac{2^{-k'\vee-k}}{2^{-k'\vee-k}+\|\sigma(x)-x_Q\|}\right)^{\theta\varepsilon} \left|2^{k\alpha}\textbf{E}_kf(x_Q)\right|^\theta
	\\\\ \ds
	= \sum_{\ell=0}^{\infty}\sum_{Q\in A_\ell} 2^{-|k'-k|\varepsilon\theta} 2^{(k'-k)\alpha\theta}  \frac{\omega(Q)^{\theta-1}}{\omega(B(\sigma(x),2^{-k'\vee-k}+\|\sigma(x)-x_Q\|))^{\theta-1}} \left(\frac{2^{-k'\vee-k}}{2^{-k'\vee-k}+\|\sigma(x)-x_Q\|}\right)^{\theta\varepsilon}
	\\\\ \ds~~~~~ 
	\times~\frac{1}{\omega(B(\sigma(x),2^{-k'\vee-k}+\|\sigma(x)-x_Q\|))} \omega(Q) \left|2^{k\alpha}\textbf{E}_kf(x_Q)\right|^\theta
	\\\\ \ds
	\leq C~2^{-|k'-k|\varepsilon\theta} 2^{(k'-k)\alpha\theta}  2^{(-k-M_0-(-k'\vee-k))N(\theta-1)} \sum_{\ell=0}^{\infty} \frac{\omega(B(\sigma(x),2^{-k'\vee-k}))^{\theta-1}}{\omega(B(\sigma(x),2^{\ell-1+(-k'\vee-k)}))^{\theta-1}} \left(\frac{1}{2^{\ell-1}}\right)^{\theta\varepsilon} 
	\\\\ \ds~~~~~
	\times~\frac{1}{\omega(B(\sigma(x),2^{\ell-1+(-k'\vee-k)}))} \int_{\|\sigma(x)-z\|\leq2\times2^{\ell+(-k\vee-k')}} \sum_{Q\in A_\ell}\left|2^{k\alpha}\textbf{E}_kf(x_Q)\right|^\theta \chi_Q(z) d\omega(z)
	\\\\ \ds
	\leq C~ 2^{-M_0 N(\theta-1)} 2^{-|k'-k|\varepsilon\theta} 2^{(-k-(-k'\vee-k))N(\theta-1)} 2^{(k'-k)\alpha\theta} 
	\\\\ \ds~~~~~
	\times~\sum_{\ell=0}^{\infty} \frac{1}{2^{(\ell-1)(N(\theta-1)+\theta\varepsilon)}}  \M\left(\sum_{Q\in Q_{cw}^k}\left|2^{k\alpha}\textbf{E}_kf(x_Q)\right|^\theta\chi_Q\right)(\sigma(\cdot))
	\\\\ \ds
	\leq C~2^{-M_0 N(\theta-1)} 2^{-|k'-k|\varepsilon\theta} 2^{(-k-(-k'\vee-k))N(\theta-1)+(k'-k)\alpha\theta} \M\left(\sum_{Q\in Q_{cw}^k}\left|2^{k\alpha}\textbf{E}_kf(x_Q)\right|^\theta\chi_Q\right)(\sigma(\cdot))
\end{array}
\eeq
where $\M$ denote the Hardy-Littlewood maximal operator on $(\R^n,\|\cdot\|,\omega)$.
Therefore,
\bel{4.12}
\begin{array}{lr}\ds
	~~~~\left|\sum_{Q\in Q_{cw}^k} 2^{k'\alpha} \omega(Q) (D_{k'} \widetilde{E}_k)(x_{Q'},x_Q) \textbf{E}_k(f)(x_Q)\right|\chi_{Q'}(x)
	\\\\ \ds
	\leq C~\sum_{\sigma\in G} 2^{-M_0 N(1-\frac{1}{\theta})}2^{-|k'-k|\varepsilon} 2^{(-k-(-k'\vee-k))N(1-\frac{1}{\theta})+(k'-k)\alpha} \left\{\M\left(\sum_{Q\in Q_{cw}^k}\left|2^{k\alpha}\textbf{E}_kf(x_Q)\right|^\theta\chi_Q\right)(\sigma(\cdot))\right\}^{\frac{1}{\theta}} \chi_{Q'}(x).
\end{array}
\eeq

For $|\alpha|<\varepsilon$, $\max\left\{\frac{N}{N+\varepsilon},\frac{N}{N+\alpha+\varepsilon}\right\}<\theta<\min\left\{p, q, 1\right\}$, it is obvious that
\bel{4.13}
\sup_{k'} \sum_{k\in\Z} 2^{-|k-k'|\varepsilon} 2^{(-k-(-k'\vee-k))N(1-\frac{1}{\theta})+(k'-k)\alpha}<\infty.
\eeq

For $1<q<\infty$, by H\"older's inequality, we have
\bel{q>1}
\begin{array}{lr}\ds
	~~~~\left|\sum_{k\in \Z} \sum_{Q\in Q_{cw}^k} 2^{k'\alpha} \omega(Q) (D_{k'} \widetilde{E}_k)(x_{Q'},x_Q) \textbf{E}_k(f)(x_Q)\right|^q \chi_{Q'}(x)
	\\\\ \ds
	\leq C~\left\{\sum_{\sigma\in G} \sum_{k\in \Z} 2^{-M_0  N(1-\frac{1}{\theta})}2^{-|k-k'|\varepsilon}  2^{(-k-(-k'\vee-k))N(1-\frac{1}{\theta})+(k'-k)\alpha} \right.
	\\\\ \ds~~~~
	\left.\times~\left\{\M\left(\sum_{Q\in Q_{cw}^k}\left|2^{k\alpha}\textbf{E}_kf(x_Q)\right|^\theta\chi_Q\right)(\sigma(\cdot))\right\}^{\frac{1}{\theta}}\right\}^q \chi_{Q'}(x)
	\\\\ \ds
	\leq C~2^{-M_0 N(1-\frac{1}{\theta})q}\left\{\sum_{\sigma\in G} \sum_{k\in \Z} 2^{-|k-k'|\varepsilon}  2^{(-k-(-k'\vee-k))N(1-\frac{1}{\theta})+(k'-k)\alpha}\right\}^{\frac{q}{q'}}
	\\\\ \ds~~~~
	\times~\left\{\sum_{\sigma\in G} \sum_{k\in \Z} 2^{-|k-k'|\varepsilon}  2^{(-k-(-k'\vee-k))N(1-\frac{1}{\theta})+(k'-k)\alpha} \right.
	\\\\ \ds~~~~
	\left.\times~\left\{\M\left(\sum_{Q\in Q_{cw}^k}\left|2^{k\alpha}\textbf{E}_kf(x_Q)\right|^\theta\chi_Q\right)(\sigma(\cdot))\right\}^{\frac{q}{\theta}}\right\} \chi_{Q'}(x)
	\\\\ \ds
	\leq C~2^{-M_0 N(1-\frac{1}{\theta})q}\sum_{\sigma\in G} \sum_{k\in \Z} 2^{-|k-k'|\varepsilon}  2^{(-k-(-k'\vee-k))N(1-\frac{1}{\theta})+(k'-k)\alpha} 
	\\\\ \ds~~~~
	\times~\left\{\M\left(\sum_{Q\in Q_{cw}^k}\left|2^{k\alpha}\textbf{E}_kf(x_Q)\right|^\theta\chi_Q\right)(\sigma(\cdot))\right\}^{\frac{q}{\theta}} \chi_{Q'}(x).
\end{array}
\eeq

For $\max\left\{\frac{N}{N+\varepsilon},\frac{N}{N+\alpha+\varepsilon}\right\}<q\leq1$, by $q$-inequality we have
\bel{q<1}
\begin{array}{lr}\ds
	~~~~~\left|\sum_{k\in \Z} \sum_{Q\in Q_{cw}^k} 2^{k'\alpha} (D_{k'} \widetilde{E}_k)(x_{Q'},x_Q) \textbf{E}_k(f)(x_Q)\right|^q \chi_{Q'}(x)
	\\\\ \ds
	\leq C~2^{-M_0 N(1-\frac{1}{\theta})q} \sum_{\sigma\in G} \sum_{k\in \Z} 2^{-|k-k'|\varepsilon q} 2^{(-k-(-k'\vee-k))N(1-\frac{1}{\theta})q+(k'-k)\alpha q}
	\\\\ \ds~~~~
	\times~\left\{\M\left(\sum_{Q\in Q_{cw}^k}\left|2^{k\alpha}\textbf{E}_kf(x_Q)\right|^\theta\chi_Q\right)(\sigma(\cdot))\right\}^{\frac{q}{\theta}} \chi_{Q'}(x).   
\end{array}
\eeq

For $|\alpha|<\varepsilon$, $\max\left\{\frac{N}{N+\varepsilon},\frac{N}{N+\alpha+\varepsilon}\right\}<\theta<\min\left\{p,q,1\right\}$, it is obvious that
\bel{}
\sup_k \sum_{k'\in \Z} 2^{-|k-k'|\varepsilon} 2^{(-k-(-k'\vee-k))N(1-\frac{1}{\theta})+(k'-k)\alpha}<\infty,
\eeq
then 
\bel{}
\begin{array}{lr}\ds
	~~~~\left\{ \sum_{k'\in\Z} \sum_{Q'\in Q_d^{k'}} \left| \sum_{k\in\Z} \sum_{Q\in Q_{cw}^k} 2^{k'\alpha} \omega(Q) (D_{k'} \widetilde{E}_k)(x_{Q'},x_Q) \textbf{E}_k(f)(x_Q) \right|^q \chi_{Q'}(x) \right\}^{\frac{1}{q}}
	\\\\ \ds
	\leq C~2^{-M_0 N(1-\frac{1}{\theta})} \left\{\sum_{\sigma\in G} \sum_{k\in\Z} \left\{\M\left(\sum_{Q\in Q_{cw}^k}\left|2^{k\alpha}\textbf{E}_kf(x_Q)\right|^\theta\chi_Q\right)(\sigma(\cdot))\right\}^{\frac{q}{\theta}}\right\}^{\frac{1}{q}}.
\end{array}
\eeq
By using Fefferman-Stein vectored maximal function inequality (see \cite{F-S}) with $\theta<\min\left\{p,q,1\right\}$, we have
\bel{}
\begin{array}{lr}\ds~~~~
	\left\| \left\{ \sum_{k'\in \Z} \sum_{Q'\in Q_d^{k'}} \left| \sum_{k\in\Z} \sum_{Q\in Q_{cw}^k} 2^{k'\alpha} \omega(Q) (D_{k'} \widetilde{E}_k)(x_{Q'},x_Q) \textbf{E}_k(f)(x_Q) \right|^q \chi_{Q'}(\cdot) \right\}^{\frac{1}{q}} \right\|_p
	\\\\ \ds
	\leq C~ 2^{-M_0 N(1-\frac{1}{\theta})}\sum_{\sigma\in G}\left\|\left\{ \sum_{k\in \Z} \sum_{Q\in Q_{cw}^k} \left(2^{k\alpha} \left|\textbf{E}_k(f)(x_Q)\right|\right)^q \chi_{Q}(\sigma(\cdot)) \right\}^{\frac{1}{q}}  \right\|_p
	\\\\ \ds
	\leq C~2^{-M_0 N(1-\frac{1}{\theta})}\left\| \left\{ \sum_{k\in\Z} \sum_{Q\in Q_{cw}^k} \left(2^{k\alpha} \left|\textbf{E}_k(f)(x_Q)\right|\right)^q \chi_{Q}(\cdot) \right\}^{\frac{1}{q}}  \right\|_p
\end{array}		
\eeq
where the last inequality follows from the fact that $G$ is a finite group and 
\bel{}
\int_{\R^n}f(\sigma(x))d\omega(x)=\int_{\R^n}f(x)d\omega(x).
\eeq

Thus, 
\bel{}
\|\textbf{S}^\alpha_q(f)\|_p\leq C~2^{-M_0N(1-\frac{1}{\theta})}\|\textbf{S}^\alpha_{q,cw}(f)\|_p.
\eeq

{\bfseries Step 2: } $\|\textbf{S}_{q,cw}^\alpha(\textbf{R}_1(f))\|_p \leq C \|\textbf{R}_1\|_{dcz} \|\textbf{S}_{q,cw}^\alpha(\textbf{R}_1(f))\|_p$.

The $\L^2$ boundedness of $\textbf{R}_1$ together with the discrete Calder\'on reproducing formula of $f\in\L^2$ on space of homogeneous type given in (\ref{cwfE}) yields
\bel{}
\begin{array}{lr}\ds
	\|\textbf{S}_{q,cw}^\alpha(\textbf{R}_1(f))\|_p
	=\left\| \left\{ \sum_{k'\in\Z} \sum_{Q'\in Q_{cw}^{k'}} \left(2^{k'\alpha} \left|\textbf{E}_{k'}(\textbf{R}_1(f))(x_{Q'})\right|\right)^q \chi_{Q'}(\cdot) \right\}^{\frac{1}{q}} \right\|_p
	\\\\ \ds~~~~~~~~~~~~~~~~~~~~~~~~~~
	=\left\| \left\{ \sum_{k'\in \Z} \sum_{Q'\in Q_{cw}^{k'}} \left(2^{k'\alpha} \left|\sum_{k\in\Z} \sum_{Q\in Q_{cw}^k} \omega(Q) (E_{k'}R_1\widetilde{E}_k)(x_{Q'},x_Q) \textbf{E}_{k}(f)(x_{Q})\right|\right)^q \chi_{Q'}(\cdot ) \right\}^{\frac{1}{q}}  \right\|_p.
\end{array}
\eeq
Applying the almost orthogonal estimate given in Lemma \ref{almost orthogonality 2} to $(E_{k'}R_1\widetilde{E}_k)(x_{Q'},x_Q)$, we obtain that for $|\alpha|<\varepsilon<1$, 
\bel{}
\begin{array}{lr}\ds
	~~~~\left|(E_{k'}R_1\widetilde{E}_k)(x_{Q'},x_Q)\right| \chi_{Q'}(x)
	\\\\ \ds
	\leq C~ \|\textbf{R}_1\|_{dcz} 2^{-|k'-k|\varepsilon} \frac{1}{V(x_{Q'},x_Q,2^{-k'\vee-k}+d(x_{Q'},x_Q))} \left(\frac{2^{-k'\vee-k}}{2^{-k'\vee-k}+d(x_{Q'},x_Q)}\right)^\varepsilon \chi_{Q'}(x)
	\\\\ \ds
	\leq C~ \|\textbf{R}_1\|_{dcz} \sum_{\sigma\in G} 2^{-|k'-k|\varepsilon} \frac{1}{V(x_{Q'},x_Q,2^{-k'\vee-k}+\|\sigma(x_{Q'})-x_Q\|)} \left(\frac{2^{-k'\vee-k}}{2^{-k'\vee-k}+\|\sigma(x_{Q'})-x_Q\|}\right)^\varepsilon \chi_{Q'}(x)
	\\\\ \ds
	\leq C~ \|\textbf{R}_1\|_{dcz} \sum_{\sigma\in G} 2^{-|k'-k|\varepsilon} \frac{1}{\omega(B(\sigma(x),2^{-k'\vee-k}+\|\sigma(x)-x_Q\|))} \left(\frac{2^{-k'\vee-k}}{2^{-k'\vee-k}+\|\sigma(x)-x_Q\|}\right)^\varepsilon \chi_{Q'}(x)
\end{array}
\eeq
where we use the fact that if $x\in Q'$, then $2^{-k\vee-k'}+\|\sigma(x_{Q'})-x_Q\|\sim2^{-k\vee-k'}+\|\sigma(x)-x_Q\|$ and $\omega(B(x_{Q'}, 2^{-k\vee-k'}+\|\sigma(x)-x_Q\|))\chi_{Q'}(x)\sim\omega(B(\sigma(x),2^{-k\vee-k'}+\|\sigma(x)-x_Q\|))\chi_{Q'}(x)$.

Similar to the estimate in \textbf{Step 1}, for $\max\left\{\frac{N}{N+\varepsilon},\frac{N}{N+\varepsilon+\alpha}\right\}<\theta<\min\left\{p, q, 1\right\}$ we have
\bel{}
\begin{array}{lr}\ds~~~~~
	\left|\sum_{Q\in Q_{cw}^k}2^{k'\alpha}\omega(Q)(E_{k'}R_1\widetilde{E}_k)(x_{Q'},x_Q)\textbf{E}_k(f)(x_Q)\right|\chi_{Q'}(x)
	\\\\ \ds
	\leq C~\|\textbf{R}_1\|_{dcz}\sum_{\sigma\in G}2^{-M_0 N(1-\frac{1}{\theta})} 2^{-|k'-k|\varepsilon} 2^{(-k-(-k\vee-k'))N(1-\frac{1}{\theta})+(k'-k)\alpha} 
	\\\\ \ds~~~~~
	\times~\left\{\M\left(\sum_{Q\in Q_{cw}^k}\left|2^{k\alpha}\textbf{E}_{k'}(f)(x_Q)\right|^\theta\chi_{Q}\right)(\sigma(\cdot))\right\}^{\frac{1}{\theta}}\chi_{Q'}(x).
\end{array}
\eeq

For $1<q<\infty$, by H\"older's inequality, we have
\bel{2.31}
\begin{array}{lr}\ds
	~~~~~\left|\sum_{k\in\Z} \sum_{Q\in Q_{cw}^k} 2^{k'\alpha}\omega(Q) (E_{k'}R_1\widetilde{E}_k)(x_{Q'},x_Q) \textbf{E}_{k}(f)(x_{Q})\right|^q \chi_{Q'}(x)
	\\\\ \ds
	\leq C~2^{-M_0 N(1-\frac{1}{\theta})q} \|\textbf{R}_1\|_{dcz}^q \sum_{\sigma\in G} \sum_{k\in\Z}  2^{-|k-k'|\varepsilon} 2^{(-k-(-k\vee-k'))N(1-\frac{1}{\theta})+(k'-k)\alpha} 
	\\\\ \ds~~~~
	\times~\left\{\M\left(\sum_{Q\in Q_{cw}^k}\left|2^{k\alpha}\textbf{E}_{k}(f)(x_Q)\right|^\theta\chi_{Q}\right)(\sigma(\cdot))\right\}^{\frac{q}{\theta}} \chi_{Q'}(x).
\end{array}
\eeq

For $\max\left\{\frac{N}{N+\varepsilon},\frac{N}{N+\alpha+\varepsilon}\right\}<q\leq1$, by $q$-inequality we have
\bel{2.32}
\begin{array}{lr}\ds
	~~~~~\left| \sum_{k\in\Z} \sum_{Q\in Q_{cw}^k} 2^{k'\alpha} \omega(Q) (E_{k'}R_1\widetilde{E}_k)(x_{Q'},x_Q) \textbf{E}_{k}(f)(x_{Q})\right|^q \chi_{Q'}(x)
	\\\\ \ds
	\leq C~2^{-M_0 N(1-\frac{1}{\theta})q} \|\textbf{R}_1\|_{dcz}^q \sum_{\sigma\in G} \sum_{k\in\Z} 2^{-|k-k'|\varepsilon q} 2^{(-k-(-k\vee-k'))N(1-\frac{1}{\theta})q+(k'-k)\alpha q} 
	\\\\ \ds~~~~
	\times~\left\{\M\left(\sum_{Q\in Q_{cw}^k}\left|2^{k\alpha}\textbf{E}_{k}(f)(x_Q)\right|^\theta\chi_{Q}\right)(\sigma(\cdot))\right\}^{\frac{q}{\theta}} \chi_{Q'}(x).
\end{array}
\eeq

Then (\ref{2.31})-(\ref{2.32}) implies that  
\bel{}
\begin{array}{lr}\ds~~~
	\left\{ \sum_{k'\in\Z} \sum_{Q'\in Q_{cw}^{k'}} \left|  \sum_{k\in\Z} \sum_{Q\in Q_{cw}^k} 2^{k'\alpha} \omega(Q) (E_{k'}R_1\widetilde{E}_k)(x_{Q'},x_Q) \textbf{E}_{k}(f)(x_{Q})\right|^q \chi_{Q'}(x)\right\}^{\frac{1}{q}}
	\\\\ \ds
	\leq C~2^{-M_0 N(1-\frac{1}{\theta})} \|\textbf{R}_1\|_{dcz} \left\{\sum_{\sigma\in G} \sum_{k\in\Z}  \left\{\M\left(\sum_{Q\in Q_{cw}^k}\left|2^{k\alpha}\textbf{E}_{k}(f)(x_Q)\right|^\theta\chi_{Q}\right)(\sigma(\cdot))\right\}^{\frac{q}{\theta}}\right\}^{\frac{1}{q}}.
\end{array}
\eeq

The Fefferman-Stein vector valued maximal function inequality with $\theta<\min\left\{p,q, 1\right\}$ yields
\bel{}
\begin{array}{lr}\ds	
	\|\textbf{S}_{q,cw}^\alpha(\textbf{R}_1(f))\|_p \leq C~2^{-M_0 N(1-\frac{1}{\theta})}\|\textbf{R}_1\|_{dcz} \left\|\left\{ \sum_{k\in\Z} \sum_{Q\in Q_{cw}^k} \left|2^{k\alpha}\textbf{E}_k(f)(x_Q)\right|^q \chi_Q(\sigma(\cdot)) \right\}^{\frac{1}{q}}\right\|_p
	\\\\ \ds~~~~~~~~~~~~~~~~~~~~~~~~~~
	\leq C~2^{-M_0 N(1-\frac{1}{\theta})} \|\textbf{R}_1\|_{dcz} \|\textbf{S}_{q,cw}^\alpha(f)\|_p.
\end{array}
\eeq

Applying the similar proof, we also have 
\bel{}
\|\textbf{S}^\alpha_{q,cw}(\textbf{R}_2(f))\|_p\leq C~2^{-M_0 N(1-\frac{1}{\theta})}\|\textbf{R}_2\|_{dcz}\|\textbf{S}^\alpha_{q,cw}(f)\|_p.
\eeq

{\bfseries Step 3: } $\|\textbf{S}^\alpha_{q,cw}(f)\|_p\leq C\|\textbf{S}^\alpha_q(f)\|_p$.

To show this estimate, the key point is to write
\bel{f}
f(x)=\textbf{T}_M(f)(x)+\textbf{R}_1(f)(x)+\textbf{R}_2(f)(x).
\eeq
Recall the estimates in \textbf{Step 2} for $\max\left\{\frac{N}{N+\varepsilon}, \frac{N}{N+\alpha+\varepsilon}\right\}\leq p<\infty$, 
\bel{cwR1}
\|\textbf{S}^{\alpha}_{q,cw}(\textbf{R}_1(f))\|_p\leq C~2^{-M_0 N(1-\frac{1}{\theta})}\|\textbf{R}_1\|_{dcz}\|\textbf{S}^{\alpha}_{q,cw}(f)\|_p,
\eeq
and
\bel{cwR2}
\|\textbf{S}^{\alpha}_{q,cw}(\textbf{R}_2(f))\|_p\leq C~2^{-M_0 N(1-\frac{1}{\theta})}\|\textbf{R}_2\|_{dcz}\|\textbf{S}^{\alpha}_{q,cw}(f)\|_p.
\eeq
Since $\|\textbf{R}_1\|_{dcz}+\|\textbf{R}_2\|_{dcz}\leq C~2^{-M\delta}$, $\delta>0$, and we choose $M$ is sufficiently larger than $M_0$, then we can set $C~2^{-M_0N(1-\frac{1}{\theta})}(\|\textbf{R}_1\|_{dcz}+\|\textbf{R}_2\|_{dcz})\leq\min\left\{\left({1\over 2}\right)^{1\over p}, {1\over 2}\right\}$, then for $\max\left\{\frac{N}{N+\varepsilon}, \frac{N}{N+\alpha+\varepsilon}\right\}<p\leq1$, by using $p$-inequality, we have
\bel{cwf}
\|\textbf{S}^\alpha_{q,cw}(f)\|_p^p\leq\|\textbf{S}^\alpha_{q,cw}(\textbf{T}_M(f)+\textbf{R}_1(f)+\textbf{R}_2(f))\|_p^p\leq \|\textbf{S}^\alpha_{q,cw}(\textbf{T}_M(f))\|_p^p+{1\over 2}\|\textbf{S}^\alpha_{q,cw}(f)\|_p^p,
\eeq
and for $1<p<\infty$, by using Minkowski's inequality we have
\bel{}
\begin{array}{lr}\ds
	\|\textbf{S}^\alpha_{q,cw}(f)\|_p\leq\|\textbf{S}^\alpha_{q,cw}(\textbf{T}_M(f))\|_p+\|\textbf{S}^\alpha_{q,cw}(\textbf{R}_1(f))\|_p+\|\textbf{S}^\alpha_{q,cw}(\textbf{R}_2(f))\|_p
	\\\\ \ds~~~~~~~~~~~~~~~~~~~
	\leq C~\|\textbf{S}^\alpha_{q,cw}(\textbf{T}_M(f))\|_p+{1\over 2}\|\textbf{S}^\alpha_{q,cw}(f)\|_p. 
\end{array}
\eeq
Hence, 
\bel{cwfcwTM}
\|\textbf{S}^\alpha_{q,cw}(f)\|_p\leq C_p \|\textbf{S}^\alpha_{q,cw}(\textbf{T}_M(f))\|_p.
\eeq

\textbf{Claim:} $\|\textbf{S}^\alpha_{q,cw}(\textbf{T}_M(f))\|_p\leq C_p \|\textbf{S}^\alpha_q(f)\|_p$.

Indeed, observing that
\bel{TM}
\textbf{T}_M(f)(x)=\sum_{k\in\Z}\sum_{Q\in Q_d^k}\omega(Q)D_k^M(x,x_Q)\textbf{D}_k(f)(x_Q),
\eeq
and
\bel{Ek'TM}
\begin{array}{lr}\ds
	\left|\textbf{E}_{k'}\textbf{T}_M(f)(x)\right|\leq \sum_{k\in\Z}\sum_{Q\in Q_d^k}\omega(Q)\left|(E_{k'}D_k^M)(x,x_Q)\right|~\left|\textbf{D}_k(f)(x_Q)\right|.
\end{array}
\eeq
Following the same proof as in \textbf{Step 1}, there exists a constant $C>0$, such that 
\bel{Ek'DkM}
\begin{array}{lr}\ds~~~~
	\left|(E_{k'}D_k^M)(x,x_Q)\right|\chi_{Q'}(x)
	\\\\ \ds
	\leq C\sum_{\sigma\in G}2^{-|k-k'|\varepsilon}{1\over \omega(B(\sigma(x),2^{-k\vee-k'}+\|\sigma(x)-x_Q\|))}\left({2^{-k\vee-k'}\over2^{-k\vee-k'}+\|\sigma(x)-x_Q\|}\right)^{\varepsilon}\chi_{Q'}(x).
\end{array}
\eeq
Therefore for $\max\left\{\frac{N}{N+\varepsilon},\frac{N}{N+\alpha+\varepsilon}\right\}<\theta<\min\{p,q,1\}$, we have
\bel{sum Q estimate}
\begin{array}{lr}\ds~~~~
	\left|\sum_{Q\in Q_d^k}2^{k'\alpha}\omega(Q)(E_{k'}D_k^M)(x,x_Q)\textbf{D}_k(f)(x_Q)\right|\chi_{Q'}(x)
	\\\\ \ds
	\leq C\sum_{\sigma\in G}2^{-M N(1-\frac{1}{\theta})}2^{-|k-k'|\varepsilon}2^{(-k-(-k'\vee-k))N(1-{1\over \theta})+(k'-k)\alpha}\left\{\M\left(\sum_{Q\in Q_d^k}\left|2^{k\alpha}\textbf{D}_k(f)(x_Q)\right|^{\theta}\right)\left(\sigma(\cdot)\right)\right\}^{1\over\theta}\chi_{Q'}(x).
\end{array}
\eeq

For $1<q<\infty, \max\{{N\over N+\varepsilon}, {N\over N+\varepsilon+\alpha}\}<\theta<\min\left\{p,q,1\right\}$, by H\"older's inequality, we obtain 
\bel{sum kQ estimate q>1}
\begin{array}{lr}\ds~~~~
	\left|\sum_{k\in\Z}\sum_{Q\in Q_d^k}2^{k'\alpha}\omega(Q)(E_{k'}D_k^M)(x,x_Q)\textbf{D}_k(f)(x_Q)\right|^q\chi_{Q'}(x)
	\\\\ \ds
	\leq C~2^{-M N(1-\frac{1}{\theta})q}\sum_{\sigma\in G}\sum_{k\in \Z}2^{-|k-k'|\varepsilon}2^{(-k-(-k'\vee-k))N(1-{1\over \theta})+(k'-k)\alpha}
	\\\\ \ds~~~~~
	\times~\left\{\M\left(\sum_{Q\in Q_d^k}\left|2^{k\alpha}\textbf{D}_k(f)(x_Q)\right|^{\theta}\right)\left(\sigma(\cdot)\right)\right\}^{q\over\theta}\chi_{Q'}(x).
\end{array}
\eeq

For $\max\left\{\frac{N}{N+\varepsilon},\frac{N}{N+\alpha+\varepsilon}\right\}<q\leq1$, by $q$-inequaliaty, we have
\bel{sum kQ estimate q<1}
\begin{array}{lr}\ds~~~~
	\left|\sum_{k\in\Z}\sum_{Q\in Q_d^k}2^{k'\alpha}\omega(Q)(E_{k'}D_k^M)(x,x_Q)\textbf{D}_k(f)(x_Q)\right|^q\chi_{Q'}(x)
	\\\\ \ds
	\leq C~2^{-M N(1-\frac{1}{\theta})q}\sum_{\sigma\in G}\sum_{k\in \Z}2^{-|k-k'|\varepsilon q}2^{(-k-(-k'\vee-k))N(1-{1\over \theta})q+(k'-k)\alpha q}
	\\\\ \ds~~~~~
	\times~\left\{\M\left(\sum_{Q\in Q_d^k}\left|2^{k\alpha}\textbf{D}_k(f)(x_Q)\right|^{\theta}\right)\left(\sigma(\cdot)\right)\right\}^{q\over\theta}\chi_{Q'}(x).
\end{array}
\eeq
Thus (\ref{sum kQ estimate q>1})-(\ref{sum kQ estimate q<1}) implies that
\bel{ScwTMf}
\begin{array}{lr}\ds~~~~
	\left\{\sum_{k'\in\Z}\sum_{Q'\in Q_{cw}^{k'}}\left|\sum_{k\in\Z}\sum_{Q\in Q_d^k}2^{k'\alpha}\omega(Q)(E_{k'}D_k^M)(x,x_Q)\textbf{D}_k(f)(x_Q)\right|^q\chi_{Q'}(x)\right\}^{1\over q}
	\\\\ \ds
	\leq C~2^{-M N(1-\frac{1}{\theta})}\left\{\sum_{\sigma\in G}\sum_{k\in \Z}\left\{\M\left(\sum_{Q\in Q_d^k}\left|2^{k\alpha}\textbf{D}_k(f)(x_Q)\right|^{\theta}\right)\left(\sigma(\cdot)\right)\right\}^{q\over\theta}\right\}^{1\over q}.
\end{array}
\eeq

The Fefferman-Stein vector valued maximal function inequality with $\theta<\min\left\{p,q,1\right\}$ yields
\bel{ScwTMf norm}
\|\textbf{S}^{\alpha}_{q,cw}(\textbf{T}_M(f))\|_{p}\leq C~2^{-MN(1-\frac{1}{\theta})}\left\|\left\{\sum_{k\in \Z}\sum_{Q\in Q^k_d}\left(2^{k\alpha}\left|\textbf{D}_k(f)(x_Q)\right|\right)^q\chi_Q(\cdot)\right\}^{1\over q}\right\|_p\leq C~\|\textbf{S}^\alpha_q(f)\|_p,
\eeq
where $M$ is a fixed constant. 

The proof of \textbf{Step 3} is complete.

Observing that $f(x)=\textbf{T}_M(f)(x)+\textbf{R}_1(f)(x)+\textbf{R}_2(f)(x)$ and applying the above estimates, we have $\|\textbf{S}^{\alpha}_q(\textbf{R}_1+\textbf{R}_2)(f)\|_p\leq C~(\|\textbf{R}_1\|_{dcz}+\|\textbf{R}_2\|_{dcz})\|\textbf{S}^{\alpha}_q(f)\|_p$, so
\bel{I-TM norm}
\|\textbf{S}^{\alpha}_q(\textbf{I}-\textbf{T}_M)(f)\|_p = C~\|\textbf{S}^{\alpha}_q(\textbf{R}_1+\textbf{R}_2)(f)\|_p \leq {1\over 2}~\|\textbf{S}^{\alpha}_q(f)\|_p.
\eeq
Similarly, for $(\textbf{T}_M)^{-1}(f)$, we also have $\|\textbf{S}^{\alpha}_q(\textbf{T}_M(f))^{-1}\|_p\leq C\|\textbf{S}^\alpha_q(f)\|_p$. If there exist a function $h\in\L^2$ and set $h=(\textbf{T}_M)^{-1}f$, we obtain
\bel{h reproducing formula}
f(x)=\textbf{T}_M(h)(x)=\sum_{k\in\Z}\sum_{Q\in Q_d^k}\omega(Q)D_{k}^M(x,x_Q)\textbf{D}_k(h)(x_Q),
\eeq
where $\|f\|_2\sim\|h\|_2$, and $\|f\|_{\dot{F}^{\alpha,q}_p}\sim\|h\|_{\dot{F}^{\alpha,q}_p}$, $|\alpha|<1$, $\max\left\{\frac{N}{N+1},\frac{N}{N+\alpha+1}\right\}<p<\infty$, $\max\left\{\frac{N}{N+1}, \frac{N}{N+\alpha+1}\right\}<q<\infty$. 

It remains to show that the series (\ref{h reproducing formula}) converges in $\L^2\cap\dot{F}^{\alpha, q}_p$. To this end, we only need to prove
\bel{converges norm}
\left\|\textbf{S}^{\alpha}_q\left(\sum_{|k|>m}\sum_{Q\in Q_d^k}\omega(Q)D_{k}^M(x,x_Q)\textbf{D}_k(h)(x_Q)\right)\right\|_p\rightarrow0,~~~\textrm{as}~~~m\rightarrow\infty.
\eeq
Repeating the same proof in \textbf{Step 1},
\bel{converges norm estimate}
\left\|\textbf{S}^{\alpha}_q\left(\sum_{|k|>m}\sum_{Q\in Q_d^k}\omega(Q)D_{k}^M(x,x_Q)\textbf{D}_k(h)(x_Q)\right)\right\|\leq C\left\|\left\{\sum_{|k|>m}\sum_{Q\in Q_d^k}\left(2^{k\alpha}\left|\textbf{D}_k(h)(x_Q)\right|\right)^q\chi_Q(x)\right\}^{1\over q}\right\|
\eeq
where by the fact $\|h\|_{\dot{F}^{\alpha,q}_p}\sim\|f\|_{\dot{F}^{\alpha,q}_p}$, the last term tends to $0$ as $m\rightarrow\infty$.

The proof of Theorem \ref{Thm 3} is complete.

\section{Duality estimates}
\setcounter{equation}{0}

The duality estimates will be a key idea in developing the Dunkl-Triebel-Lizorkin space theory. Now we show Theorem \ref{dual prop}. 

\textbf{Proof of Theorem \ref{dual prop} (A)}:

Applying the weak-type discrete Calder\'on-type reproducing formula given in Theorem \ref{Thm 3} for $f\in \L^2\cap\dot{F}^{\alpha,q}_p$, we write
\bel{}
f(x)=\sum_{k\in\Z}\sum_{Q\in Q_d^k}\omega(Q)D_k^M(x,x_Q)\textbf{D}_k(h)(x_Q),
\eeq
where $\|h\|_2\sim\|f\|_2$, and $\|\textbf{S}^\alpha_q(h)\|_p\sim\|\textbf{S}^\alpha_q(f)\|_p$.

For $1<q<\infty$, by H\"older's inequality
\bel{}
\begin{array}{lr}\ds
	\left|\left<f,g\right>\right| =\left|\left<\sum_{k\in\Z}\sum_{Q\in Q_d^k}\omega(Q)D_k^M(\cdot,x_Q)\textbf{D}_k(h)(x_Q), g(\cdot)\right>\right|
	\\\\ \ds~~~~~~~~~~~~
	\leq\int_{\R^n} \sum_{k\in\Z} \sum_{Q\in Q_d^k} 2^{k\alpha} \left|\textbf{D}_k(h)(x_Q)\right| 2^{-k\alpha} \left|\textbf{D}_k^M(g)(x_Q)\right| \chi_Q(x) d\omega(x)
	\\\\ \ds~~~~~~~~~~~~
	\leq C \int_{\R^n} \left\{\sum_{k\in\Z} \sum_{Q\in Q_d^k} \left(2^{k\alpha}\left|\textbf{D}_k(h)(x_Q)\right|\right)^q \chi_Q(x)\right\}^{\frac{1}{q}} 
	\\\\ \ds~~~~~~~~~~~~~~~~~
	\times\left\{\sum_{k\in\Z} \sum_{Q\in Q_d^k} \left(2^{-k\alpha} \left|\textbf{D}_k^M(g)(x_Q)\right|\right)^{q'} \chi_Q(x)\right\}^{\frac{1}{q'}} d\omega(x)
	\\\\ \ds~~~~~~~~~~~~
	\leq C\left\{\int_{\R^n}\left|\textbf{S}^{\alpha}_q(h)(x)\right|^pd\omega(x)\right\}^{1\over p}\left\{\int_{\R^n}\left|\textbf{S}^{-\alpha}_{q'}(g)(x)\right|^{p'}d\omega(x)\right\}^{1\over p'}
	\\\\ \ds~~~~~~~~~~~~
	\leq C~\|f\|_{\dot{F}^{\alpha,q}_p}~\|g\|_{\dot{F}^{-\alpha,q'}_{p'}}.
\end{array}
\eeq

\textbf{Proof of Theorem \ref{dual prop} (B):}

Set $\Omega_\ell=\left\{x\in\R^n:\textbf{S}_q^\alpha(h)(x)>2^\ell\right\}$, and 
$B_\ell=\left\{Q:\omega(Q\cap\Omega_\ell)>\frac{1}{2}\omega(Q)~\textrm{and}~{\omega(Q\cap\Omega_{\ell+1})}\leq\frac{1}{2}\omega(Q)\right\}$ where $Q$ are dyadic cubes. Denote $B^\ast_\ell:=\left\{Q^\ast_\ell\right\}$ as the maximal dyadic cubes in $B_\ell$ for $\ell\in\Z$. 

\textbf{Claim: }the Calder\'on reproducing formula can be rewrite as:
\bel{}
f(x)=\sum_{\ell\in\Z}\sum_{Q^\ast_\ell\in B^\ast_\ell}\sum_{Q\subset Q^\ast_{\ell}}\omega(Q)D_k^M(x,x_Q)\textbf{D}_k(h)(x_Q).
\eeq
In order to prove the above claim, we only need to show that if the dyadic cube $Q\notin B_\ell$ for all $\ell\in\Z$, then
\bel{0}
\omega(Q)D_k^M(x,x_Q)\textbf{D}_k(h)(x_Q)=0.
\eeq

Observe that by the stopping time argument, each dyadic cube $Q$ can be in one and only one $B_\ell$, that is, if $Q$ belongs to both $B_\ell$ and $B_{\ell'}$, then $\ell=\ell'$. We now assume that $\omega(Q)\neq0$. Otherwise, the equality (\ref{0}) holds obviously. Note that $\omega(\Omega_\ell)<2^{-2\ell}\|\textbf{S}^\alpha_q(h)\|_2^2\rightarrow0$ as $\ell\rightarrow+\infty$. As a consequence, if $Q\notin B_\ell$ for all $\ell\in\Z$, then $\omega(Q\cap\Omega_\ell)\leq{1\over2}\omega(Q)$ for all $\ell\in\Z$ since, otherwise, there exists an $\ell_0\in\Z$, such that $\omega(Q\cap\Omega_{\ell_0})>{1\over 2}\omega(Q)$. However, $\omega(Q\cap\Omega_\ell)\rightarrow0$ as $\ell\rightarrow+\infty$ and $\{\omega(Q\cap\Omega_\ell)\}_\ell$ is a decreasing sequence. So there must be a critical index $\ell_1$ such that $\omega(Q\cap\Omega_{\ell_1})>{1\over 2}\omega(Q)$ and $\omega(Q\cap\Omega_{\ell_1+1})\leq{1\over 2}\omega(Q)$, that is $Q\in B_{\ell_1}$. This is contradict to the fact that $Q$ is not in $B_\ell$ for all $\ell\in\Z$.

Since $\omega(Q\cap\Omega_\ell)\leq{1\over 2}\omega(Q)$ for all $\ell\in\Z$, then $\omega(Q\cap\Omega_\ell^c)\geq{1\over 2}\omega(Q)$ for all $\ell\in\Z$. Set $K=\left\{{x\in\R^n},\textbf{S}^\alpha_q(h)(x)=0\right\}$. Note that $\cap_{\ell\in \Z}\Omega_\ell^c=\cap_{\ell\in\Z}\left\{x\in\R^n:\textbf{S}^\alpha_q(h)(x)\leq2^\ell\right\}=K$. Thus
\bel{}
\omega(Q\cap K)=\lim_{\ell\rightarrow-\infty}\omega(Q\cap\Omega_\ell^c)\geq{1\over 2}\omega(Q)>0
\eeq
for all $x\in K$, $0=\textbf{S}^\alpha_q(h)(x)=\left\{\sum_{k\in\Z}\sum_{Q\in Q_d^k}\left(2^{k\alpha}\left|\textbf{D}_k(h)(x_Q)\right|\right)^q\chi_Q(x)\right\}^{1\over q}$, then $\left|\textbf{D}_k(h)(x_Q)\right|=0$. Hence the claim is proved.

Then by H\"{o}lder's inequality and $p$-inequality with $p\leq1$, we have
\bel{omega}
\begin{array}{lr}\ds
	\left|\left<f,g\right>\right|= \left| \sum_{\ell\in\Z} \sum_{Q^\ast_\ell\in B^\ast_\ell} \sum_{Q\subset Q^\ast_\ell} \omega(Q)^{\frac{1}{q}} 2^{k\alpha} \textbf{D}_k(h)(x_Q)   \omega(Q)^{\frac{1}{q'}} 2^{-k\alpha} \textbf{D}_k^M(g)(x_Q) \right|
	\\\\ \ds~~~~~~~~~~~~
	\leq \sum_{\ell\in \Z} \sum_{Q^\ast_\ell\in B^\ast_\ell} \left( \sum_{Q\subset Q^\ast_\ell} \omega(Q) \left|2^{k\alpha} \textbf{D}_k(h)(x_Q)\right|^q \right)^{\frac{1}{q}}   \left( \sum_{Q\subset Q^\ast_\ell} \omega(Q) \left|2^{-k\alpha} \textbf{D}_k^M(g)(x_Q)\right|^{q'} \right)^{\frac{1}{q'}}
	\\\\ \ds~~~~~~~~~~~~
	\leq \left\{ \sum_{\ell\in\Z} \sum_{Q^\ast_\ell\in B_\ell^\ast} \left( \sum_{Q\subset Q^\ast_\ell} \omega(Q) \left|2^{k\alpha} \textbf{D}_k(h)(x_Q)\right|^q \right)^{\frac{p}{q}}   \left( \sum_{Q\subset Q^\ast_\ell} \omega(Q) \left|2^{-k\alpha} \textbf{D}_k^M(g)(x_Q)\right|^{q'} \right)^{\frac{p}{q'}} \right\}^{\frac{1}{p}}
	\\\\ \ds~~~~~~~~~~~~
	= \left\{\sum_{\ell\in\Z} \sum_{Q^\ast_\ell\in B^\ast_\ell} \left( \sum_{Q\subset Q^\ast_\ell} \omega(Q) \left|2^{k\alpha} \textbf{D}_k(h)(x_Q)\right|^q \right)^{\frac{p}{q}} \right.
	\\\\ \ds~~~~~~~~~~~~~~~~
	\left.\times~\omega(Q^\ast_\ell)^{1-\frac{p}{q}}  \left( \frac{1}{\omega(Q^\ast_\ell)^{\frac{q'}{p}-\frac{q'}{q}}} \sum_{Q\subset Q^\ast_\ell} \omega(Q) \left|2^{-k\alpha} \textbf{D}_k^M(g)(x_Q)\right|^{q'} \right)^{\frac{p}{q'}} \right\}^{\frac{1}{p}}
	\\\\ \ds~~~~~~~~~~~~
	\leq\left\{ \sum_{\ell\in\Z} \sum_{Q^\ast_\ell\in B^\ast_\ell} \left( \sum_{Q\subset Q^\ast_\ell} \omega(Q) \left|2^{k\alpha} \textbf{D}_k(h)(x_Q)\right|^q \right)^{\frac{p}{q}} \omega(Q^\ast_\ell)^{1-\frac{p}{q}} \right\}^{1\over p} \|g\|_{CMO^{-\alpha,q'}_{p}}
	\\\\ \ds~~~~~~~~~~~~
	\leq \left\{\sum_{\ell\in\Z}\left\{\sum_{Q_\ell^\ast\in B^\ast_\ell}\sum_{Q\subset Q_\ell^\ast}\omega(Q)\left|2^{k\alpha}\textbf{D}_k(h)(x_Q)\right|^q\right\}^{\frac{p}{q}}\left\{\sum_{Q_\ell^\ast\in B^\ast_\ell}\omega(Q^\ast_\ell)\right\}^{1-\frac{p}{q}}\right\}^\frac{1}{p} \|g\|_{CMO^{-\alpha,q'}_{p}}.
\end{array}
\eeq
To estimate the term $\left\{\sum_{Q_\ell^\ast\in B^\ast_\ell}\omega(Q^\ast_\ell)\right\}^{1-\frac{p}{q}}$. Set $\widetilde{\Omega}_\ell=\left\{x\in\R^n:\M(\chi_{\Omega_\ell})(x)>\frac{1}{2}\right\}$, where $\M$ is the Hardy-Littlewood maximal function on $\R^n$ with the measure $d\omega$ and $\chi_{\Omega_\ell}(x)$ is the indicate function of $\Omega_\ell$. It is easily to see that if $Q\in B_\ell$, then $Q\subset\widetilde{\Omega}_\ell$. Since all $Q^\ast_\ell$ are disjoint, thus
\bel{}
\left\{\sum_{Q_\ell^\ast\in B^\ast_\ell}\omega(Q_\ell^\ast)\right\}^{1-\frac{p}{q}}\leq C~\omega(\widetilde{\Omega}_\ell)^{1-\frac{p}{q}}\leq C~\omega(\Omega_\ell)^{1-{p\over q}},
\eeq 
where the first inequality follows from the fact that $\cup_{Q^\ast_\ell\in B^\ast_\ell}Q^\ast_\ell\subset\widetilde{\Omega}_\ell$ and $\sum_{Q^\ast_\ell\in B^\ast_\ell}\omega(Q^\ast_\ell)\leq\omega(\widetilde{\Omega}_\ell)$, and by the $\L^2$-boundedness of the Hardy-Littlewood maximal function, the last inequality follows from the estimate $\omega(\widetilde{\Omega}_\ell)\leq C\omega(\Omega_\ell)$.

We claim that
\bel{claim 1}
\begin{array}{lr}\ds
	\sum_{Q^\ast_\ell\in B^\ast_\ell}\sum_{Q\subset Q^\ast_\ell}\omega(Q)\left|2^{k\alpha}\textbf{D}_k(h)(x_Q)\right|^q \leq C~2^{q\ell}\omega(\Omega_\ell).	
\end{array}
\eeq
Under this claim (\ref{claim 1}), we get 
\bel{}
\begin{array}{lr}\ds
	\left|\left<f,g\right>\right|\leq C\left(\sum_{\ell\in\Z}\left(2^{q\ell}\omega(\Omega_\ell)\right)^{\frac{p}{q}}\omega(\Omega_{\ell})^{1-\frac{p}{q}}\right)^{\frac{1}{p}}\|g\|_{CMO^{-\alpha,q'}_p}
	\\\\ \ds~~~~~~~~~~~~
	\leq C\left(\sum_{\ell\in\Z}2^{p\ell}\omega(\Omega_\ell)\right)^\frac{1}{p}\|g\|_{CMO^{-\alpha,q'}_p}
	\\\\ \ds~~~~~~~~~~~~
	\leq C \|\textbf{S}_q^\alpha(h)\|_p\|g\|_{CMO^{-\alpha,q'}_p}
	\\\\ \ds~~~~~~~~~~~~
	\leq C \|f\|_{\dot{F}^{\alpha, q}_{p}}\|g\|_{CMO^{-\alpha,q'}_p}.
\end{array}
\eeq
It remains to show the claim (\ref{claim 1}). 
In order to do that, we begin with the following estimate
\bel{}
\begin{array}{lr}\ds
	\int_{\widetilde{\Omega}_\ell\setminus\Omega_{\ell+1}}\textbf{S}_q^\alpha(h)(x)^qd\omega(x)\leq C~2^{q\ell} \omega(\widetilde{\Omega}_\ell)\leq C~ 2^{q\ell} \omega({\Omega_\ell}).
\end{array}
\eeq
Note that
\bel{}
\begin{array}{lr}\ds
	\int_{\widetilde{\Omega}_\ell\setminus\Omega_{\ell+1}}\left|\textbf{S}_q^\alpha(h)(x)\right|^qd\omega(x)\geq \sum_{Q\in B_\ell}\left(2^{k\alpha}\left|\textbf{D}_k(h)(x_Q)\right|\right)^q\omega\left((\widetilde{\Omega}_\ell\setminus\Omega_{\ell+1})\cap Q\right).
\end{array}
\eeq

Since for each $Q\in B_\ell$ implies $Q\subseteq\widetilde{\Omega}_\ell$ and $\Omega_{\ell+1}\subset\Omega_{\ell}$. Thus
\bel{}
\omega\left((\widetilde{\Omega}_\ell\setminus\Omega_{\ell+1})\cap Q\right)=\omega(Q)-\omega(\Omega_{\ell+1}\cap Q)\geq\frac{1}{2}\omega(Q).
\eeq
Therefore,
\bel{}
\int_{\widetilde{\Omega}_\ell\setminus\Omega_{\ell+1}}\left|\textbf{S}_q^\alpha(h)(x)\right|^qd\omega(x)\geq C\sum_{Q\in B_\ell}\omega(Q)\left(2^{k\alpha}\left|\textbf{D}_k(h)(x_Q)\right|\right)^q.
\eeq
This implies the claim (\ref{claim 1}). The prove of Theorem \ref{dual prop} (\ref{0<p<1,1<q<infty}) is complete.	

\textbf{Proof of Theorem \ref{dual prop} (C)}:

Recall Definition \ref{Fap infty}, for $g\in\L^2(\R^n,\omega)$, the norm of $g\in\dot{F}^{\alpha,\infty}_p(\R^n,\omega)$ is defined by
\bel{}
\|f\|_{\dot{F}^{\alpha,\infty}_{p}}=\left\|{\sup_{{k\in\Z},{Q\in Q_d^k}}}2^{k\alpha}\left|\textbf{D}_k^M(f)(x_Q)\right|\chi_{Q}(x)\right\|_{\L^p(\R^n,\omega)}.
\eeq 

By H\"older's inequality, we have
\bel{dual 3}
\begin{array}{lr}\ds
	\left|\left<f,g\right>\right|=\left|\sum_{k\in\Z}\sum_{Q\in Q_d^k}\omega(Q)\textbf{D}_k^M(g)(x_Q)\textbf{D}_k(h)(x_Q)\right|
	\\\\ \ds~~~~~~~~~~~~
	\leq C \int_{\R^n}\sum_{k\in\Z}\sum_{Q\in Q_d^k}2^{k\alpha}\left|\textbf{D}_k(h)(x_Q)\right|\chi_Q(x)\textbf{S}^{-\alpha}_{\infty}(g)(x)d\omega(x)
	\\\\ \ds~~~~~~~~~~~~
	\leq C~\|\textbf{S}^\alpha_1(h)\|_p~\|\textbf{S}^{-\alpha}_{\infty}(g)\|_{p'}
	\\\\ \ds~~~~~~~~~~~~
	\leq C~ \|f\|_{\dot{F}^{\alpha,1}_p}~\|g\|_{\dot{F}^{-\alpha,\infty}_{p'}}.
\end{array}
\eeq

To get Theorem \ref{dual prop} (C), it suffices to show $\|f\|_{\dot{F}^{\alpha,1}_p}\leq C~\|f\|_{\dot{F}^{\alpha,q}_p}$.

According to the $q$-inequality with $q\leq 1$, we find
\bel{3}
\begin{array}{lr}\ds
	\|f\|_{\dot{F}^{\alpha,1}_p}= C\left\{\int_{\R^n}\left(\sum_{k\in\Z}\sum_{Q\in Q_d^k}2^{k\alpha}\left|\textbf{D}_k(h)(x_Q)\right|\chi_Q(x)\right)^p d\omega(x)\right\}^{1\over p}
	\\\\ \ds~~~~~~~~~~~~
	\leq C\left\{\int_{\R^n}\left|\sum_{k\in\Z}\sum_{Q\in Q_d^k}\left(2^{k\alpha}\left|\textbf{D}_k(h)(x_Q)\right|\right)^q\chi_Q(x)\right|^{p\over q} d\omega(x)\right\}^{1\over p}
	\\\\ \ds~~~~~~~~~~~~
	\leq C~\|\textbf{S}^\alpha_q(h)\|_p
	
	\leq C~\|f\|_{\dot{F}^{\alpha,q}_p}.
\end{array}
\eeq
As a result,
\bel{}
\left|\left<f,g\right>\right|\leq C \|f\|_{\dot{F}^{\alpha,q}_p}~\|g\|_{\dot{F}^{-\alpha,\infty}_{p'}}.
\eeq

\textbf{Proof of Theorem \ref{dual prop} (D)}:

Recall Definition \ref{Fa infty infty}, for $g\in\L^2(\R^n,\omega)$, the norm of $g\in\dot{F}^{\alpha,\infty}_{\infty}(\R^n,\omega)$ is defined by
\bel{}
\|g\|_{\dot{F}^{\alpha,\infty}_{\infty}}={\sup_{{k\in\Z},{Q\in Q_d^k}}} 2^{k\alpha}\left|\textbf{D}_k^M(g)(x_Q)\right|.
\eeq

Denote $\beta=-\alpha+N(\frac{1}{p}-1)$, then we have
\bel{dual 4}
\begin{array}{lr}\ds
	\left|\left<f,g\right>\right|=\left|\sum_{k\in\Z}\sum_{Q\in Q_d^k}\omega(Q)\textbf{D}_k^M(g)(x_Q)\textbf{D}_k(h)(x_Q)\right|
	\\\\ \ds~~~~~~~~~~~~
	\leq C \int_{\R^n}\sum_{k\in\Z}\sum_{Q\in Q_d^k}2^{-k\beta}\left|\textbf{D}_k(h)(x_Q)\right|\chi_Q(x)d\omega(x) \sup_{{k\in\Z},{Q\in Q_d^k}}2^{k\beta}\left|\textbf{D}_k^M(g)(x_Q)\right|
	\\\\ \ds~~~~~~~~~~~~
	= C~ \|S_1^{-\beta}(h)\|_{1}~\|g\|_{\dot{F}^{\beta,\infty}_\infty}
	\\\\ \ds~~~~~~~~~~~~
	\leq C~ \|f\|_{\dot{F}^{\alpha-N({1\over p}-1),1}_1}~\|g\|_{\dot{F}^{-\alpha+N({1\over p}-1),\infty}_\infty}.
\end{array}
\eeq

Recall that $\omega(B(x,r))\geq C~r^N$ deduced by (\ref{sim}). Since $p\leq1$, and $\ell(Q)=2^{-k-M}$, where $M$ is a fixed constant, then $2^{-kN({1\over p}-1)}\leq C{\omega(B(x,2^{-k}))}^{{1\over p}-1}\sim C~\omega(Q)^{{1\over p}-1}$. Similar to the proof of Theorem \ref{dual prop} (C),it suffices to show $\|f\|_{\dot{F}^{\alpha-N({1\over p}-1),1}_1}\leq C~\|f\|_{\dot{F}^{\alpha,q}_p}$. By Minkowski's inequality, and $p$-inequality with $p\leq1$, we find
\bel{5.21}
\begin{array}{lr}\ds
	\|f\|_{\dot{F}^{\alpha-N(\frac{1}{p}-1),1}_1}= \int_{\R^n}\sum_{k\in\Z}\sum_{Q\in Q_d^k}2^{k(\alpha-N(\frac{1}{p}-1))}\left|\textbf{D}_k(f)(x_Q)\right|\chi_{Q}(x)d\omega(x)
	\\\\ \ds~~~~~~~~~~~~~~~~~~~~~
	\leq C  \sum_{k\in\Z}\sum_{Q\in Q_d^k}2^{k\alpha}\omega(Q)^{\frac{1}{p}-1}\left|\textbf{D}_k(f)(x_Q)\right|\omega(Q)
	\\\\ \ds~~~~~~~~~~~~~~~~~~~~~
	=C  \sum_{k\in\Z}\sum_{Q\in Q_d^k}2^{k\alpha}\left|\textbf{D}_k(f)(x_Q)\right|\omega(Q)^{\frac{1}{p}}
\end{array}
\eeq

If $\frac{p}{q}\geq1$, applying Minkowski's inequality and $p$-inequality again, the right-hand side of $(\ref{5.21})$ is further controlled by   
\bel{}
\begin{array}{lr}\ds~~~~
	\left(\sum_{k\in\Z}\sum_{Q\in Q_d^k}2^{k\alpha p}\left|\textbf{D}_k(f)(x_Q)\right|^p\omega(Q)\right)^{\frac{1}{p}}
	\\\\ \ds
	\leq\left\{\int_{\R^n}\sum_{k\in\Z}\sum_{Q\in Q_d^k}\left(2^{k\alpha q}\left|\textbf{D}_k(f)(x_Q)\right|^q\chi_Q(x)\right)^{\frac{p}{q}}d\omega(x)\right\}^{1\over p}
	\\\\ \ds
	\leq \left\{\int_{\R^n} \left\{\sum_{k\in\Z}\sum_{Q\in Q_d^k}2^{k\alpha q}\left|\textbf{D}_k(f)(x_Q)\right|^q\chi_Q(x)\right\}^{\frac{p}{q}}d\omega(x)\right\}^{\frac{1}{p}}
	\\\\ \ds
	= C~\|\textbf{S}^{\alpha}_q(f)\|_p = C~\|f\|_{\dot{F}^{\alpha,q}_p}.
\end{array}
\eeq

If $\frac{p}{q}<1$, applying $q$-inequality and the Minkowski's inequality again, the estimates of right-hand side of $(\ref{5.21})$ are following  
\bel{}
\begin{array}{lr}\ds~~~~
	\sum_{k\in\Z}\sum_{Q\in Q_d^k}2^{k\alpha}\left|\textbf{D}_k(f)(x_Q)\right|\omega(Q)^{\frac{1}{p}}
	\\\\ \ds
	=\sum_{k\in\Z}\sum_{Q\in Q_d^k}\left(2^{k\alpha p}\left|\textbf{D}_k(f)(x_Q)\right|^p\omega(Q)\right)^{\frac{1}{p}}
	\\\\ \ds
	\leq\left\{\sum_{k\in\Z}\sum_{Q\in Q^k_d}\left(2^{k\alpha p}\left|\textbf{D}_k(f)(x_Q)\right|^p\omega(Q)\right)^{q\over p}\right\}^{1\over q}
	\\\\ \ds
	=\left\{\left\{\sum_{k\in\Z}\sum_{Q\in Q^k_d}\left(\int_{\R^n}2^{k\alpha p}\left|\textbf{D}_k(f)(x_Q)\right|^p\chi_Q(x)d\omega(x)\right)^{q\over p}\right\}^{p\over q}\right\}^{1\over p}
	\\\\ \ds
	\leq\left\{\int_{\R^n}\left\{\sum_{k\in\Z}\sum_{Q\in Q^k_d}2^{k\alpha q}\left|\textbf{D}_k(f)(x_Q)\right|^q\chi_Q(x)\right\}^{p\over q}d\omega(x)\right\}^{1\over p}
	\\\\ \ds
	= \|\textbf{S}^{\alpha}_q(f)\|_p = \|f\|_{\dot{F}^{\alpha,q}_p}.
\end{array}
\eeq

The proof of Theorem \ref{dual prop} (\ref{0<p,q<1}) is complete.

The Theorem \ref{dual prop} (\ref{1<p,q<infty}) indicates that if $\{f_n\}_{n=1}^\infty$ is a sequence in $\L^2$ with $\|\textbf{S}^\alpha_q(f_n-f_m)\|_p\rightarrow0$ as $n,m\rightarrow\infty$. Then for each $g\in\L^2$ with $\|g\|_{\dot{F}^{-\alpha,q'}_{p'}}<\infty$, $1<p<\infty$, $1<q<\infty$, we have $\lim\limits_{n,m\rightarrow\infty}\left<f_n-f_m,g\right>=0$. Therefore, there exists $f$, as a distribution on $\L^2\cap {\dot{F}^{-\alpha,q'}_{p'}}$, such that for each $g\in\L^2$ with $\|g\|_{\dot{F}^{-\alpha,q'}_{p'}}<\infty$,
\bel{}
\left<f,g\right>=\lim_{n\rightarrow\infty}\left<f_n,g\right>.
\eeq
Other ranges of $p$, $q$ stated above have the same results.

\section{Dunkl-Triebel-Lizorkin space}
\setcounter{equation}{0}

In this section, we define the Dunkl-Triebel-Lizorkin spaces, and show the spaces are complete. Before introducing the Dunkl-Triebel-Lizorkin space, we need the following 
\begin{lemma}\label{Lemma 4.1}
	\quad 
	\begin{enumerate}[(A)]
		\item \label{1<p,q<infty}
		For $1<p<\infty$, $1<q<\infty$. Then $D_k(\cdot, y)$ is in $\L^2\cap\dot{F}^{-\alpha,q'}_{p'}$ for any fixed $k$ and $y\in\R^n$. Moreover
		\bel{}
		\|\textbf{S}^{-\alpha}_{q'}\left(D_k(\cdot, y)\right)\|_{p'}\leq C,
		\eeq
		where the constant $C$ only depends on $k$.
		\item \label{0<p<1,1<q<infty}
		For $\max\left\{\frac{N}{N+1},\frac{N}{N+\alpha+1}\right\}<p\leq1$, $1<q<\infty$. Then $D_k(\cdot,y)$ is in $\L^2\cap CMO_p^{-\alpha,q'}$ for any fixed $k$ and $y\in\R^n$. Moreover
		\bel{}
		\sup_P\left( \frac{1}{\omega(P)^{\frac{q'}{p}-\frac{q'}{q}}}\sum_{Q\subset P}\omega(Q)\left( 2^{-j\alpha}\left|\textbf{D}_j(D_k(\cdot,y))(x_Q)\right|\right)^{q'}\right)^{\frac{1}{q'}} \leq C
		\eeq
		where both $P$ and $Q$ are dyadic cubes on $\R^n$ and the constant $C$ which depends on $k$ but is independent of $y$.
		\item \label{1<p<infty,0<q<1}
		For $1<p<\infty$, $\max\left\{{N\over N+1}, {N\over N+\alpha+1}\right\}<q\leq 1$. Then $D_k(\cdot, y)$ is in $\L^2\cap\dot{F}^{-\alpha, \infty}_{p'}$ for any fixed $k$ and $y\in\R^n$. Moreover
		\bel{}
		\|\textbf{S}^{-\alpha}_\infty(D_k(\cdot,y))\|_{p'}=\left\|\sup_{{j\in\Z},{Q\in Q_d^j}}2^{-j\alpha}\left|\textbf{D}_j(D_k(\cdot,y))(x_Q)\right|\chi_Q(x)\right\|_{p'}<C.
		\eeq
		where the constant $C$ only depends on $k$.
		\item \label{0<p,q<1}
		For $\max\left\{{N\over N+1}, {N\over N+\alpha+1}\right\}<p\leq1$, $\max\left\{{N\over N+1}, {N\over N+\alpha+1}\right\}<q\leq1$. Then $D_k(\cdot, y)$ is in $\L^2\cap\dot{F}^{-\alpha+N({1\over p}-1),\infty}_{\infty}$ for any fixed $k$ and $y\in\R^n$. Moreover
		\bel{}
		\sup_{{j\in\Z},{Q\in Q_d^j}}2^{(-\alpha+N(\frac{1}{p}-1))j}\left|\textbf{D}_j(D_k(\cdot,y))(x_Q)\right|\leq C.
		\eeq		
	\end{enumerate}		
\end{lemma}

\textbf{Proof of Lemma \ref{Lemma 4.1} (\ref{1<p,q<infty})}:

According to Lemma \ref{almost orthogonality 1}, we choose $|\alpha|<\varepsilon<1$ such that 
\bel{}
\left|\textbf{D}_j(D_k(\cdot, y))(x_Q)\right|\leq C~2^{-|j-k|\varepsilon}{1\over V(x_Q, y, 2^{-j\vee-k}+d(x_Q, y))}\left({2^{-j\vee-k}\over 2^{-j\vee-k}+\|x_Q-y\|}\right)^{\varepsilon}.
\eeq
Then,
\bel{}
\begin{array}{lr}\ds~~~~~
	\|\textbf{S}^{-\alpha}_{q'}\left(D_k(\cdot, y)\right)\|_{p'}
	\\\\ \ds
	=\left\{\int_{R^n}\left\{\sum_{j\in\Z}\sum_{Q\in Q_d^j}\left(2^{-j\alpha}\left|\textbf{D}_j\left(D_k(\cdot, y)\right)(x_Q)\right|\chi_{Q}(x)\right)^{q'}\right\}^{{p'}\over{q'}}d\omega(x)\right\}^{1\over p'}
	\\\\ \ds
	\leq C~\left\{\int_{R^n}\left\{\sum_{j\in\Z}\left(\sum_{Q\in Q_d^j}2^{-j\alpha}2^{-|j-k|\varepsilon}{1\over \omega(B(x,2^{-j\vee-k}+d(x,y)))}\left({2^{-j\vee-k}\over 2^{-j\vee-k}+\|x-y\|}\right)^{\varepsilon}\chi_{Q}(x)\right)^{q'}\right\}^{{p'}\over{q'}}d\omega(x)\right\}^{1\over p'}
\end{array}
\eeq
where the inequality also follows from the fact that $q'>1$ and $x\in Q$.

If ${{p'}\over{q'}}\geq 1$, by Minkowski's inequality, we have
\bel{}
\begin{array}{lr}\ds~~~~~
	\|\textbf{S}^{-\alpha}_{q'}(D_k(\cdot,y))\|_{p'}
	\\\\ \ds
	\leq C~\left\{\sum_{j\in\Z}\left\{\int_{\R^n}\left(2^{-j\alpha} 2^{-|j-k|\varepsilon}{1\over \omega(B(x, 2^{-j\vee-k}+d(x,y)))}\left({2^{-j\vee-k}\over 2^{-j\vee-k}+\|x-y\|}\right)^{\varepsilon}\right)^{p'}d\omega(x)\right\}^{{q'}\over {p'}}\right\}^{1\over q'}
	\\\\ \ds
	\leq C \sum_{j\in\Z}2^{-j\alpha} 2^{-|j-k|\varepsilon}\left\{\int_{\R^n}\left({1\over \omega\left(x, 2^{-j\vee-k}+d(x,y)\right)}\left({2^{-j\vee-k}\over 2^{-j\vee-k}+\|x-y\|}\right)^{\varepsilon}\right)^{p'}d\omega(x)\right\}^{1\over p'}.
\end{array} 
\eeq

If ${{p'}\over{q'}}<1$, applying Minkowski's inequality again, and ${1\over p'}$-inequality with $p'>1$, we have
\bel{}
\begin{array}{lr}\ds~~~~~
	\|\textbf{S}^{-\alpha}_{q'}(D_k(\cdot, y))\|_{p'}
	\\\\ \ds
	\leq C~\left\{\int_{\R^n}\sum_{j\in\Z}\left(2^{-j\alpha}2^{-|j-k|\varepsilon}{1\over \omega(B(x, 2^{-j\vee-k}+d(x,y)))}\left({2^{-j\vee-k}\over 2^{-j\vee-k}+\|x-y\|}\right)^{\varepsilon}\right)^{p'}d\omega(x)\right\}^{{1}\over {p'}}
	\\\\ \ds
	\leq C~\left\{\sum_{j\in\Z}\int_{\R^n}\left(2^{-j\alpha}2^{-|j-k|\varepsilon}{1\over \omega(B(x, 2^{-j\vee-k}+d(x,y)))}\left({2^{-j\vee-k}\over 2^{-j\vee-k}+\|x-y\|}\right)^{\varepsilon}\right)^{p'}d\omega(x)\right\}^{1\over p'}
	\\\\ \ds
	\leq C \sum_{j\in\Z}2^{-j\alpha}2^{-|j-k|\varepsilon}\left\{\int_{\R^n}\left({1\over \omega(B(x, 2^{-j\vee-k}+d(x,y)))}\left({2^{-j\vee-k}\over 2^{-j\vee-k}+\|x-y\|}\right)^{\varepsilon}\right)^{p'}d\omega(x)\right\}^{1\over p'}.
\end{array}
\eeq
Thus,
\bel{1}
\begin{array}{lr}\ds~~~~
	\|\textbf{S}^{-\alpha}_{q'}\left(D_k(\cdot, y)\right)\|_{p'}
	\\\\ \ds
	\leq C\sum_{j\in\Z}2^{-j\alpha}2^{-|j-k|\varepsilon}
	\left\{\int_{\R^n}\left({1\over \omega(B(x, 2^{-j\vee-k}+d(x,y)))}\left({2^{-j\vee-k}\over 2^{-j\vee-k}+\|x-y\|}\right)^{\varepsilon}\right)^{p'}d\omega(x)\right\}^{1\over p'}.
\end{array}
\eeq

For $j\geq k$, 
\bel{2.86}
\begin{array}{lr}\ds~~~~
	\int_{\R^n}\left({1\over \omega(B(x, 2^{-j\vee-k}+d(x,y)))}\left({2^{-j\vee-k}\over 2^{-j\vee-k}+\|x-y\|}\right)^{\varepsilon}\right)^{p'}d\omega(x)
	\\\\ \ds
	\leq C \int_{d(x,y)<2^{-k}}\left({1\over \omega(B(x, 2^{-j\vee-k}+d(x,y)))}\left({2^{-j\vee-k}\over 2^{-j\vee-k}+\|x-y\|}\right)^{\varepsilon}\right)^{p'}d\omega(x)
	\\\\ \ds~~~~~
	+\sum_{m=1}^{\infty}\int_{2^{-k+(m-1)}\leq d(x,y)<2^{-k+m}}\left({1\over \omega(B(x, 2^{-j\vee-k}+d(x,y)))}\left({2^{-j\vee-k}\over 2^{-j\vee-k}+\|x-y\|}\right)^{\varepsilon}\right)^{p'}d\omega(x).
\end{array}
\eeq
Note that
\bel{}
\begin{array}{lr}\ds~~~~~
	\int_{d(x,y)<2^{-k}}\left({1\over \omega(B(x, 2^{-j\vee-k}+d(x,y)))}\left({2^{-j\vee-k}\over 2^{-j\vee-k}+\|x-y\|}\right)^{\varepsilon}\right)^{p'}d\omega(x)
	\\\\ \ds
	\leq C \sum_{\sigma\in G}\int_{\|\sigma(y)-x\|<2^{-k}}{1\over \omega(B(x, 2^{-k}))^{p'}}d\omega(x)
	\\\\ \ds
	\leq C \sum_{\sigma\in G}\int_{\|\sigma(y)-x\|<2^{-k}}{1\over \omega(B(\sigma(y), 2^{-k}))^{p'}}d\omega(x)
	\\\\ \ds
	\leq C~\omega(B(\sigma(y), 2^{-k}))^{1-p'}
	\\\\ \ds
	\leq C~2^{-kN(1-p')}.
\end{array}
\eeq
And for $m\geq1$,
\bel{}
\begin{array}{lr}\ds~~~~~
	\int_{2^{-k+(m-1)}\leq d(x,y)<2^{-k+m}}\left({1\over \omega(B(x, 2^{-j\vee-k}+d(x,y)))}\left({2^{-j\vee-k}\over 2^{-j\vee-k}+\|x-y\|}\right)^{\varepsilon}\right)^{p'}d\omega(x)
	\\\\ \ds
	\leq C~2^{-m\varepsilon}\sum_{\sigma\in G}\int_{\|\sigma(y)-x\|<2^{-k+m}}{1\over \omega(B(\sigma(y), 2^{-k+m}))^{p'}}d\omega(x)
	\\\\ \ds
	\leq C~2^{-m\left(\varepsilon+(p'-1)N\right)}2^{-kN(1-p')}.
\end{array}
\eeq
Thus, for $j\geq k$,  
\bel{}
\begin{array}{lr}\ds
	~~~~\int_{\R^n}\left({1\over \omega(B(x, 2^{-j\vee-k}+d(x,y)))}\left({2^{-j\vee-k}\over 2^{-j\vee-k}+\|x-y\|}\right)^{\varepsilon}\right)^{p'}d\omega(x)\leq C~2^{-kN(1-p')}\leq C_k,
\end{array}
\eeq
where $C_k$ is a constant depending on $k$.

Similarly, for $j<k$,
\bel{}
\begin{array}{lr}\ds
	~~~~\int_{\R^n}\left({1\over \omega(B(x, 2^{-j\vee-k}+d(x,y)))}\left({2^{-j\vee-k}\over 2^{-j\vee-k}+\|x-y\|}\right)^{\varepsilon}\right)^{p'}d\omega(x)\leq C~2^{-jN(1-p')}.
\end{array}
\eeq
Together with (\ref{1}) and $|\alpha|<\varepsilon$, we have 
\bel{}
\begin{array}{lr}\ds
	\|\textbf{S}^{-\alpha}_{q'}\left(D_k(\cdot, y)\right)\|_{p'}\leq C_k\sum_{j:j\geq k}2^{-j\alpha}2^{-(j-k)\varepsilon}+\sum_{j:j< k}2^{-j\alpha}2^{(j-k)\varepsilon}2^{-jN({1\over p'}-1)}\leq C_k.
\end{array}
\eeq

\textbf{Proof of Lemma \ref{Lemma 4.1} (\ref{0<p<1,1<q<infty})}:

Fix a dyadic cube $P$ with the side length $2^{-M-j_0}$ and the center $x_P$. For $Q\in Q_d^j$, applying Lemma \ref{almost orthogonality 1}, we have
\bel{2.92}
\begin{array}{lr}\ds
	\left| \textbf{D}_j\left( D_k(\cdot,y) \right) (x_Q) \right|\leq C~2^{-|j-k|\varepsilon} \frac{1}{V(x_Q,y,2^{-j\vee-k}+d(x_Q,y))} \left( \frac{2^{-j\vee-k}}{2^{-j\vee-k}+\|x_Q-y\|} \right)^{\varepsilon}
	\\\\ \ds~~~~~~~~~~~~~~~~~~~~~~~~~~~~~~~
	\leq C~2^{-|j-k|\varepsilon} \frac{1}{\omega(B(x_Q,2^{-j\vee-k}))} 
	\\\\ \ds~~~~~~~~~~~~~~~~~~~~~~~~~~~~~~~
	\leq C~2^{-|j-k|\varepsilon} \frac{1}{\omega(B(x_Q,2^{-k}))}.
\end{array}
\eeq
We have the estimate
\bel{}
\begin{array}{lr}\ds
	~~~~\frac{1}{\omega(P)^{\frac{q'}{p}-\frac{q'}{q}}}\sum_{Q\subset P}\omega(Q)\left( 2^{-j\alpha}\left|\textbf{D}_j(D_k(\cdot,y))(x_Q)\right|\right)^{q'}
	\\\\ \ds
	= \frac{1}{\omega(P)^{\frac{q'}{p}-\frac{q'}{q}}}\sum_{j=j_0}^{\infty}\sum_{\{Q\in Q_d^j:Q\subset P\}}\omega(Q)\left( 2^{-j\alpha}\left|\textbf{D}_j(D_k(\cdot,y))(x_Q)\right|\right)^{q'}
	\\\\ \ds
	\leq C \frac{1}{\omega(P)^{\frac{q'}{p}-\frac{q'}{q}}}\sum_{j=j_0}^{\infty}\sum_{\{Q\in Q_d^j:Q\subset P\}}\omega(Q)2^{-j\alpha q'}2^{-|j-k|\varepsilon q'}\frac{1}{\omega(B(x_Q,2^{-k}))^{q'}}
	\\\\ \ds
	\leq C\sup_{x\in P}\frac{1}{\omega(B(x,2^{-k}))^{q'}} \frac{1}{\omega(P)^{\frac{q'}{p}-\frac{q'}{q}}}\sum_{j=j_0}^{\infty}\sum_{\{Q\in Q_d^j:Q\subset P\}}\omega(Q)2^{-(j\alpha+|j-k|\varepsilon)q'}
	\\\\ \ds
	\leq C\sup_{x\in P}\frac{1}{\omega(B(x,2^{-k}))^{q'}} \frac{1}{\omega(P)^{\frac{q'}{p}-\frac{q'}{q}}}\sum_{j=j_0}^{\infty}\omega(P)2^{-(j\alpha+|j-k|\varepsilon)q'}
	\\\\ \ds
	\leq C\sup_{x\in P}\frac{1}{\omega(B(x,2^{-k}))^{q'}} \frac{1}{\omega(P)^{\frac{q'}{p}-\frac{q'}{q}-1}}\sum_{j=j_0}^{\infty}2^{-(j\alpha+|j-k|\varepsilon)q'}.
\end{array}
\eeq

\begin{enumerate}[(i)]
	\item 
	For $j_0\geq k$, the doubling property of the measure $\omega$ implies
	\bel{}
	\omega(B(x_P,2^{-k}))\leq C~2^{(-k+j_0)N}\omega(B(x_P ,2^{-j_0}))\sim C~2^{(-k+j_0)N}\omega(P).
	\eeq 
	Thus, for $|\alpha|<\varepsilon$,
	\bel{}
	\begin{array}{lr}\ds
		~~~~\frac{1}{\omega(P)^{\frac{q'}{p}-\frac{q'}{q}}}\sum_{Q\subset P}\omega(Q)\left( 2^{-j\alpha}\left|\textbf{D}_j(D_k(\cdot,y))(x_Q)\right|\right)^{q'}
		\\\\ \ds
		\leq C\sup_{x\in P}\frac{1}{\omega(B(x,2^{-k}))^{q'}} \frac{1}{\omega(B(x_P,2^{-k}))^{\frac{q'}{p}-\frac{q'}{q}-1}}2^{(-k+j_0)N(\frac{q'}{p}-\frac{q'}{q}-1)}\sum_{j=j_0}^{\infty}2^{-((j-k)\varepsilon+j\alpha)q'}
		\\\\ \ds
		\leq C_k\sup_{x\in P}\frac{1}{\omega(B(x,2^{-k}))^{\frac{q'}{p}}}.
	\end{array}
	\eeq
	\item 
	For $j_0<k$, since 
	\bel{}
	\omega(B(x_Q,2^{-k}))\leq\omega(B(x_Q,2^{-j_0}))\sim\omega(P), 
	\eeq 
	we have
	\bel{}
	\begin{array}{lr}\ds~~~~
		\frac{1}{\omega(P)^{\frac{q'}{p}-\frac{q'}{q}}}\sum_{Q\subset P}\omega(Q)\left( 2^{-j\alpha}\left|\textbf{D}_j(D_k(\cdot,y))(x_Q)\right|\right)^{q'}
		\\\\ \ds
		\leq C\sup_{x\in P}\frac{1}{\omega(B(x,2^{-k}))^{q'}} \frac{1}{\omega(P)^{\frac{q'}{p}-\frac{q'}{q}-1}}\sum_{j=j_0}^{\infty}2^{-(j\alpha+|j-k|\varepsilon)q'}
		\\\\ \ds
		\leq C_k\sup_{x\in P}\frac{1}{\omega(B(x,2^{-k}))^{\frac{q'}{p}}}.
	\end{array}
	\eeq
\end{enumerate}

Taking the supremum over all dyadic cubes $P$ and using $\inf_{x\in\R^n}\omega(B(x,1))\geq C$, we obtain
\bel{}
\sup_P\left( \frac{1}{\omega(P)^{\frac{q'}{p}-\frac{q'}{q}}}\sum_{Q\subset P}\omega(Q)\left( 2^{-j\alpha}\left|\textbf{D}_j(D_k(\cdot,x_Q))(y)\right|\right)^{q'}\right)^{\frac{1}{q'}} \leq C_k
\eeq
and the proof of Lemma \ref{Lemma 4.1} (\ref{0<p<1,1<q<infty}) is complete.

\textbf{Proof of Lemma \ref{Lemma 4.1} (\ref{1<p<infty,0<q<1}):}
\bel{}
\begin{array}{lr}\ds~~~~
	\|\textbf{S}^{-\alpha}_\infty(D_k(\cdot,y))\|_{p'}
	\\\\ \ds
	=\left\{\int_{\R^n}\left|\sup_{{j\in\Z},{Q\in Q^j_d}}2^{-j\alpha}\left|\textbf{D}_j\left(D_k(\cdot,y)\right)(x_Q)\right|\chi_Q(x)\right|^{p'}d\omega(x)\right\}^{1\over{p'}}
	\\\\ \ds
	\leq C~\left\{\int_{\R^n}\left|\sup_{{j\in\Z},{Q\in Q^j_d}}2^{-j\alpha}2^{-|j-k|\varepsilon}{1\over V(x_Q,y, 2^{-j\vee-k}+d(x_Q,y))}\left({2^{-j\vee-k}\over 2^{-j\vee-k}+\|x_Q-y\|}\right)^{\varepsilon}\chi_Q(x)\right|^{p'}d\omega(x)\right\}^{1\over{p'}}
	\\\\ \ds
	\leq C~\left\{\sup_{{j\in\Z},{Q\in Q^j_d}}2^{-j\alpha p'}2^{-|j-k|\varepsilon p'}\int_Q\left({1\over \omega(B(x, 2^{-j\vee-k}+d(x,y)))}\left({2^{-j\vee-k}\over2^{-j\vee-k}+\|x-y\|}\right)^{\varepsilon}\right)^{p'}d\omega(x)\right\}^{1\over p'}
	\\\\ \ds
	\leq C\sup_{j\in\Z}2^{-j\alpha}2^{-|j-k|\varepsilon}\left\{\int_{\R^n}\left({1\over \omega(B(x, 2^{-j\vee-k}+d(x,y)))}\left({2^{-j\vee-k}\over2^{-j\vee-k}+\|x-y\|}\right)^{\varepsilon}\right)^{p'}d\omega(x)\right\}^{1\over p'}
\end{array}
\eeq
Similar to the proof of Lemma \ref{Lemma 4.1} (\ref{1<p,q<infty}), we have
\bel{}
\|\textbf{S}^{-\alpha}_\infty(D_k(\cdot,y))\|_{p'}\leq C_k\sup_{j: j\geq k}2^{-j\alpha}2^{-(j-k)\varepsilon}+\sup_{j:j<k}2^{-j\alpha}2^{(j-k)\varepsilon}2^{-jN({1\over p'}-1)}\leq C_k.
\eeq

\textbf{Proof of Lemma \ref{Lemma 4.1} (\ref{0<p,q<1}):}

Recall (\ref{2.92}), we have
\bel{}
\begin{array}{lr}\ds
	\left|\textbf{D}_j(D_k(\cdot,y))(x_Q)\right|\leq C~2^{-|j-k|\varepsilon} \frac{1}{\omega(B(x_Q,2^{-k}))}.
\end{array}
\eeq

\begin{enumerate}[(i)]
	\item 
	For $j\geq k$, since $p>\frac{N}{N+\alpha+\varepsilon}$, then 
	\bel{}
	\sup_{j:j\geq k;~Q\in Q^j_d}  2^{(-\alpha+N(\frac{1}{p}-1))j}\left|\textbf{D}_j(D_k(\cdot,y))(x_Q)\right| \leq C\sup_{j:j\geq k;~Q\in Q^j_d} 2^{(-\alpha+N(\frac{1}{p}-1))j} 2^{-(j-k)\varepsilon} \frac{1}{\omega(B(x_Q,2^{-k}))} \leq C_k,
	\eeq
	
	\item 
	For $j<k$, since $p\leq1$ and $|\alpha|<\varepsilon$, then
	\bel{}
	\sup_{j:j<k;~Q\in Q^j_d}  2^{(-\alpha+N(\frac{1}{p}-1))j}\left|\textbf{D}_j(D_k(\cdot,y))(x_Q)\right| \leq C\sup_{j:j<k;~Q\in Q^j_d}2^{(-\alpha+N(\frac{1}{p}-1))j}2^{(j-k)\varepsilon} \frac{1}{\omega(B(x_Q,2^{-k}))} \leq C_k.
	\eeq
	
\end{enumerate}

The proof of Lemma \ref{Lemma 4.1} is complete.

For $1<p<\infty$, $1<q<\infty$, we denote $\L^2\cap\dot{F}^{-\alpha,q'}_{p'}$, as the subspace of $f\in\L^2$, with the norm $\|f\|_{\dot{F}^{-\alpha,q'}_{p'}}<\infty$. 
Based on the above Lemma \ref{Lemma 4.1}, if $f\in \left(\L^2\cap\dot{F}^{-\alpha,q'}_{p'}\right)'$, then $\textbf{D}_k(f)(x)$ is well defined since for each fixed $x$, $D_k(x,y)\in\L^2\cap\dot{F}^{-\alpha,q'}_{p'}$. Other ranges of $p$, $q$ have similar results.

Theorem \ref{thm} describes an important property for each distribution $f$. More precisely, it establishes the weak-type discrete Calder\'on reproducing formula in the distribution sense. Now we prove Theorem \ref{thm}.

\textbf{Proof Theorem \ref{thm} (\ref{1<p,q<infty})}

By the Theorem \ref{dual prop}, there exists $f\in \left(\L^2\cap\dot{F}^{-\alpha, q'}_{p'}\right)'$ such that for each $g\in \L^2\cap\dot{F}^{-\alpha, q'}_{p'}$,
\bel{lim dual 1}
\left<f,g\right>=\lim_{n\rightarrow\infty}\left<f_n, g\right>. 
\eeq
Observing that 
\bel{lim S}
\|\textbf{S}^\alpha_q(f-f_n)\|_p=\|\textbf{S}^{\alpha}_q(\lim_{m\rightarrow\infty}(f_m-f_n))\|_p\leq\liminf_{m\rightarrow\infty}\|\textbf{S}^{\alpha}_q(f_m-f_n)\|_p,
\eeq     
hence $\|\textbf{S}^\alpha_q(f-f_n)\|_p\rightarrow0$ as $n\rightarrow \infty$. This implies that
\bel{lim norm}
\|f\|_{\dot{F}^{\alpha,q}_p}=\|\textbf{S}^\alpha_q(f)\|_p=\lim_{n\rightarrow\infty}\|\textbf{S}^\alpha_q(f_n)\|_p=\lim_{n\rightarrow\infty}\|f_n\|_{\dot{F}^{\alpha,q}_p}<\infty.
\eeq
Applying Theorem \ref{Thm 3}, for each $f_n$ there exists an $h_n$ such that $\|f_n\|_2\sim\|h_n\|_2$ and $\|f_n\|_{\dot{F}^{\alpha,q}_p}\sim\|h_n\|_{\dot{F}^{\alpha,q}_p}$. Thus by Theorem \ref{dual prop}, there exists $h\in\left(\L^2\cap\dot{F}^{-\alpha,q'}_{p'}\right)'$ such that for each $g\in\L^2\cap\dot{F}^{-\alpha,q'}_{p'}$, we have
\bel{lim dual h}
\left<h,g\right>=\lim_{n\rightarrow\infty}\left<h_n, g\right>.
\eeq  
Therefore, 
\bel{}
\|\textbf{S}^\alpha_q(h_n-h_m)\|_p\to0,
\eeq
and 
\bel{}
\|h\|_{\dot{F}^{\alpha,q}_p}=\|\textbf{S}^\alpha_q(h)\|_p=\lim_{n\rightarrow\infty}\|\textbf{S}^\alpha_q(h_n)\|_p\sim\lim_{n\rightarrow\infty}\|\textbf{S}^\alpha_q(f_n)\|_p=\|\textbf{S}^\alpha_q(f)\|_p=\|f\|_{\dot{F}^{\alpha,q}_p}.
\eeq        
For each $g\in \L^2\cap\dot{F}^{-\alpha, q'}_{p'}$, we know that
\bel{dual in prop 4.1}
\left|\sum_{k\in\Z}\sum_{Q\in Q_d^k}\omega(Q)\textbf{D}_k^M(g)(x_Q)\textbf{D}_k(h)(x_Q)\right|\leq C~\|f\|_{\dot{F}^{\alpha,q}_p}\|g\|_{\dot{F}^{-\alpha,q'}_{p'}},
\eeq
which implies that the series $\sum_{k\in\Z}\sum_{Q\in Q_d^k}\omega(Q)D_k^M(x,x_Q)\textbf{D}_k(h)(x_Q)$ is a distribution in $\left(\L^2\cap\dot{F}^{-\alpha, q'}_{p'}\right)'$. Moreover, by the reproducing formula of $f_n$ in Theorem \ref{Thm 3}, for each $g\in \L^2\cap\dot{F}^{-\alpha, q'}_{p'}$,
\bel{distribution reproducing formula a}
\left<f,g\right>=\lim_{n\rightarrow\infty}\left<f_n,g\right>=\lim_{n\rightarrow\infty}\left<\sum_{k\in\Z}\sum_{Q\in Q_d^k}\omega(Q)D_k^M(\cdot,x_Q)\textbf{D}_k(h_n)(x_Q),g(\cdot)\right>,
\eeq
where $\|f_n\|_2\sim\|h_n\|_2$ and $\|f_n\|_{\dot{F}^{\alpha,q}_p}\sim\|h_n\|_{\dot{F}^{\alpha,q}_p}$.

By the same proof of Theorem \ref{dual prop}, we have
\bel{hn-h}
\left|~\left<\sum_{k\in\Z}\sum_{Q\in Q_d^k}\omega(Q)D_k^M(\cdot,x_Q)\textbf{D}_k(h-h_n)(x_Q), g(\cdot)\right>~\right|\leq C~\|h_n-h\|_{\dot{F}^{\alpha,q}_p}~\|g\|_{\dot{F}^{-\alpha, q'}_{p'}}.
\eeq
Since $\|h_n-h\|_{\dot{F}^{\alpha,q}_p}\rightarrow0$, as $n\rightarrow\infty$, we have
\bel{fg}
\left<f, g\right>=\lim_{n\rightarrow\infty}\left<f_n,g\right>=\left< \sum_{k\in\Z}\sum_{Q\in Q_d^k}\omega(Q)D_k^M(\cdot,x_Q)\textbf{D}_k(h)(x_Q), g(\cdot) \right>.
\eeq

The proof of Theorem \ref{thm} (\ref{0<p<1,1<q<infty}), (\ref{1<p<infty,0<q<1}), (\ref{0<p,q<1}) is similar, so we omit the details. 

The Theorem \ref{thm} (\ref{1<p,q<infty}) indicates that one can consider $\L^2\cap\dot{F}^{-\alpha,q'}_{p'}$, $|\alpha|<1$, $1<p<\infty$, $1<q<\infty$, the subspace of $f\in\L^2$ with the norm $\|f\|_{\dot{F}^{-\alpha,q'}_{p'}}<\infty$, as the test function space and $\left(\L^2\cap\dot{F}^{-\alpha,q'}_{p'}\right)'$ as the distribution space. The Dunkl-Triebel-Lizorkin space is defined by Definition \ref{definition}. We remark that in the Definition \ref{definition}, the series $\sum_{k\in\Z}\sum_{Q\in Q_d^k}\omega(Q)\lambda_QD_k^M(x,x_Q)$ with $\left\|\left\{\sum_{k\in\Z}\sum_{Q\in Q_d^k}\left(2^{k\alpha}\left|\lambda_Q\right|\right)^q\chi_Q\right\}^{1\over q}\right\|_p<\infty$ defines a distribution in $\left(\L^2\cap \dot{F}^{-\alpha,q'}_{p'}\right)'$. Indeed, applying the proof of Theorem \ref{dual prop} for each $g\in\L^2\cap\dot{F}^{-\alpha,q'}_{p'}$,
\bel{}
\left|\sum_{k\in\Z}\sum_{Q\in Q_d^k}\omega(Q)\lambda_Q\textbf{D}_k^M(g)(x_Q)\right|\leq C~\left\|\left\{\sum_{k\in\Z}\sum_{Q\in Q_d^k}\left(2^{k\alpha}\left|\lambda_Q\right|\right)^q\chi_Q\right\}^{1\over q}\right\|_p\|g\|_{\dot{F}^{-\alpha, q'}_{p'}}.
\eeq
Other ranges of $p$, $q$ has the same results.

\textbf{Proof of Theorem \ref{closure} (\ref{1<p,q<infty})}

Suppose $f\in\dot{F}^{\alpha,q}_p(\R^n,\omega)$. Then $f\in \left(\L^2\cap\dot{F}^{-\alpha,q'}_{p'}\right)'$ and $f$ has a wavelet-type decomposition $f(x)=\sum_{k\in\Z}\sum_{Q\in Q_d^k}\omega(Q)\lambda_QD_k^M(x,x_Q)$ in $\left(\L^2\cap\dot{F}^{-\alpha,q'}_{p'}\right)'$ with its norm $\left\|\left\{\sum_{k\in\Z}\sum_{Q\in Q_d^k}\left(2^{k\alpha}\left|\lambda_Q\right|\right)^q\chi_Q(x)\right\}^{1\over q}\right\|_{p}<\infty$. Set 
\bel{fn}
f_n(x)=\sum_{|k|\leq n}\sum_{Q\in Q_d^k\atop{Q\subseteq B(0,n)}}\omega(Q)\lambda_QD_k^M(x,x_Q).
\eeq
Then $f_n\in\L^2\cap\dot{F}^{\alpha,q}_p$ and $f_n$ converges to $f$ in $\left(\L^2\cap\dot{F}^{-\alpha,q'}_{p'}\right)'$ as $n\rightarrow \infty$. 

To see that $f\in \bar{~\L^2\cap\dot{F}^{\alpha,q}_p~}$, by Theorem \ref{thm}, it suffices to show that $\|f_n-f_m\|_{\dot{F}^{\alpha, q}_p}\rightarrow0$ as $n,m\rightarrow\infty$. Indeed, if let $E_n=\left\{(k,Q):|k|\leq n,Q\in Q_d^k\subseteq B(0,n)\right\}$ and $E^{c}_{n,m}=E_n\setminus E_m$ with $n\geq m$,
\bel{fn-fm}
\begin{array}{lr}\ds
	\|f_n-f_m\|_{\dot{F}^{\alpha,q}_p}=\left\|~\left(\sum_{k'\in \Z}\sum_{Q'\in Q_d^k}\left(2^{k'\alpha}\left|\textbf{D}_{k'}^M(f_n-f_m)(x_{Q'})\right|\right)^q\chi_{Q'}(x)\right)^{1\over q}~\right\|_p
	\\\\ \ds~~~~~~~~~~~~~~~~~~~~~
	\leq \left\|~\left(\sum_{k'\in\Z}\sum_{Q'\in Q_d^k}\left(2^{k'\alpha}\left|\textbf{D}_{k'}^M\left(\sum_{E^c_{n,m}}\omega(Q)\lambda_QD_k^M(\cdot, x_Q)\right)(x_{Q'})\right|~\right)^q\chi_{Q'}(x)\right)^{1\over q}~\right\|_p
	\\\\ \ds~~~~~~~~~~~~~~~~~~~~~
	\leq C\left\|~\left\{\sum_{E^c_{n,m}}\left(2^{k\alpha}\left|\lambda_Q\right|\right)^q\chi_Q(x)\right\}^{1\over q}~\right\|_p\rightarrow0,
\end{array}
\eeq
as $n,m\to\infty$, where the last inequality follows from the same proof of \textbf{Step 1} in the Theorem \ref{Thm 3} and hence, $f\in\bar{~\L^2\cap\dot{F}^{\alpha,q}_p~}$.

Conversely, if $f\in\bar{~\L^2\cap\dot{F}^{\alpha,q}_p~}$ by Theorem \ref{dual prop}, then there exists $h\in\left(\L^2\cap\dot{F}^{-\alpha,q'}_{p'}\right)'$ with $\|\textbf{S}^\alpha_q(h)\|_p\sim\|\textbf{S}^\alpha_q(f)\|_p$ such that for each $g\in\L^2\cap\dot{F}^{-\alpha,q'}_{p'}$,
\bel{fg dual}
\left<f,g\right>=\left<\sum_{k\in\Z}\sum_{Q\in Q_d^k}\omega(Q)D_k^M(\cdot,x_Q)\textbf{D}_k(h)(x_Q), g(\cdot)\right>.
\eeq
Set $\lambda_Q=\textbf{D}_k(h)(x_Q)$ with $Q\in Q_d^k$. We obtain a wavelet-type decomposition of $f$ in $\left(\L^2\cap\dot{F}^{-\alpha,q'}_{p'}\right)'$ in the distribution sense:
\bel{f wavelet-type}
f(x)=\sum_{k\in\Z}\sum_{Q\in Q_d^k}\omega(Q)\lambda_QD_k^M(x,x_Q).
\eeq
Hence, $f\in\dot{F}^{\alpha,q}_p(\R^n,\omega)$. Moreover
\bel{f norm}
\|f\|_{\dot{F}^{\alpha, q}_p}=\inf\left\{~\left\|~\left\{\sum_{k\in\Z}\sum_{Q\in Q_d^k}\left(2^{k\alpha}\left|\lambda_Q\right|\right)^q\chi_Q(x)\right\}^{1\over q}\right\|_p~\right\}\leq C~\|\textbf{S}^\alpha_q(h)\|_p\leq C~\|\textbf{S}^\alpha_q(f)\|_p.
\eeq
The proof of Theorem \ref{closure} (\ref{0<p<1,1<q<infty}), (\ref{1<p<infty,0<q<1}), (\ref{0<p,q<1}) are similar, so we omit the details.

\section{Declarations}

\quad \par 

Disclosure of potential conflicts of interest: The authors declare no conflict of interest.

Research involving Human Participants and/or Animals: Not applicable.

Informed consent: The manuscript described has not been submitted elsewhere for publication, in whole or in part.

Funding: No funding was received for conducting this study.


\newpage


\begin{thebibliography}{100}
	
	\bibitem{JA}
	J.-Ph. Anker, J. Dziuba\'nski, A. Hejna, {\it Harmonic functions, conjugate harmonic functions and the Hardy Space $H^1$ in the rational Dunkl setting}, J. Fourier Anal. Appl., {\bf 25} (2019), no. 5, 2356–2418.
	
	\bibitem{Bui}
	T. A. Bui, {\it Harmonic analysis in Dunkl settings}, ResearchGate, Jun. 2023. [Online]. Available: https://www.researchgate.net/publication/371952480. [Accessed: Aug. 1, 2024].
	
	\bibitem{B-D-L}
	T. A. Bui, X. T. Duong, F. K. Ly, {\it Maximal function characterizations for Hardy spaces on spaces of homogeneous type with finite measure and applications}, J. Funct. Anal. {\bf 278} (2020), art. 108423, 55 pp. 
	
	\bibitem{CW}
	R. R. Coifman, G. Weiss, {\it Analyse Harmonique Non-commutative sur Certains Espaces Homog\`{e}nes}, Lecture Notes in Math., Vol. {\bf 242}, Springer-Verlag, Berlin-New York, 1971, v+160pp.
	
	\bibitem{D-J-S}
	G. David., J. L. Journ\'e, S. Semmes, {\it Op\'erateurs de Calder\'on-Zygmund, fonctions para-acccr\'etives et interpolation}, Rev. Mat. Iberoam., {\bf 1} (1985), no. 4, 1-56.
	
	\bibitem{DH1}
	J. Dziuba\'nski, A. Hejna, {\it Remark on atomic decompositions for the Hardy space $H^1$ in the rational Dunkl setting}, Studia Math., {\bf 251} (2020), no. 1, 89-110.
	
	\bibitem{DH2}
	J. Dziuba\'nski, A. Hejna, {\it Remarks on Dunkl Translations of Non-radial Kernels}, Journal of Fourier Analysis and Applications, {\bf 29} (2023), no.4, Paper No. 52, 35 pp.
	
	\bibitem{F-J1}
	M. Frazier, B. Jawerth, {\it A discrete transform and decompositions of distribution spaces}, J. Funct. Anal. {\bf 93} (1990), 34–170.
	
	\bibitem{F-S}
	L. Grafakos, L. Liu, D. Yang, {\it Vector-valued singular integrals and maximal functions on spaces of homogeneous type}, Math.
	Scand., {\bf 104} (2009), no. 2, 296–310.
	
	\bibitem{Han}
	Y. Han, {\it Plancherel-P\'olya type inequality on spaces of homogeneous type and its applications}, Proc. Amer. Math.Soc., {\bf 126} (1998), no. 11, 3315-3327.
	
	\bibitem{H-L-L-W}
	Y. Han, M. Lee, J. Li, B. Wick, {\it Lipschitz and Triebel–Lizorkin spaces, commutators in Dunkl setting}, Nonlinear Anal. {\bf 237} (2023), no. 113365, 36 pp.
	
	\bibitem{H-M-Y}
	Y. Han, D. Müller, D. Yang, {\it A theory of Besov and Triebel–Lizorkin spaces on metric measure spaces modeled on Carnot–Carathéodory spaces}, Abstr. Appl. Anal. 2008, art. 893409, 250 pp.
	
	\bibitem{Han and Sawyer}
	Y. Han, E.T Sawyer, {\it Littlewood-Paley theory on spaces of homogeneous type and classical function spaces}, Mem. Amer. Math. Soc. {\bf 110} (1994), no. 530, 1-126. 
	
	\bibitem{Lizorkin1}
	P. I. Lizorkin, {\it Operators connected with fractional differentiation, and classes of differentiable functions}, Trudy Mat. Inst. Steklov. {\bf 117} (1972), no. 345, 212–243.
	
	\bibitem{Lizorkin2}
	P. I. Lizorkin, {\it Properties of functions in the spaces $\lambda_{p,\theta}^r$}, Trudy Mat. Inst. Steklov. {\bf 131} (1974), no. 247, 158–181.
	
	\bibitem{Peetre1}
	J. Peetre, {\it Remarques sur les espaces de Besov. Le cas $0<p<1$}, C. R. Acad. Sci. Paris Sér. A-B {\bf 277} (1973), 947–949.
	
	\bibitem{Peetre2}
	J. Peetre, {\it On spaces of Triebel–Lizorkin type}, Ark. Mat. {\bf 13} (1975), 123–130.
	
	\bibitem{Peetre3}
	J. Peetre, {\it New Thoughts on Besov Spaces}, Duke Univ. Math. Ser., Duke Univ. Press, Durham, 1976.
	
	\bibitem{Hp}
	C. Tan, Y. Han, Y. Han, M.-Y. Lee, J. Li, {\it Singular integral operators, $T1$ theorem, Littlewood-Paley theory and Hardy spaces in Dunkl Setting }, arXiv:2204.01886.
	
	\bibitem{TX}
	S. Thangavelua, Y. Xu, {\it Riesz transform and Riesz potentials for Dunkl transform}, Journal of Computational and Applied Mathematics, {\bf 199} (2007), 181–195.
	
	\bibitem{Triebel}
	H. Triebel, {\it Spaces of distributions of Besov type on Euclidean $n$-space. Duality, interpolation}, Ark. Mat., {\bf 11} (1973), 13–64.
	
\end{thebibliography}
\end{document}